\documentclass{article}
\usepackage{amsmath,amssymb,amsfonts,amsthm,amsopn,dsfont}
\usepackage{bm}
\usepackage[ruled,vlined,linesnumbered]{algorithm2e}
\newtheorem{prop}{Proposition}
\newtheorem{lemma}{Lemma}
\usepackage{graphicx}
\usepackage[backend=bibtex]{biblatex}

\newcommand{\hsp}{\mathcal{H}^{\Gamma}}

\newcommand{\extg}{\mathcal{E}_{_{\Gamma}}}
\newcommand{\exts}{\mathcal{E}_{_{\Sigma}}}
\newcommand{\trg}{\gamma_{_{\Gamma}}}

\newcommand{\psid}{\psi_{_D}}
\newcommand{\psis}{\psi_{_\Sigma}}
\newcommand{\hpsid}{\hat{\psi}_{_D}}
\newcommand{\hpsis}{\hat{\psi}_{_\Sigma}}
\begin{document}

	\title{A PDE-constrained optimization method for 3D-1D coupled problems with discontinuous solutions}
	
	\author{Stefano Berrone
		\footnote{stefano.berrone@polito.it}, 
		Denise Grappein
		\footnote{denise.grappein@polito.it},
		Stefano Scial\'o
		\footnote{stefano.scialo@polito.it},\\{\footnotesize Dipartimento di Scienze Matematiche, Politecnico di Torino,
			Corso Duca degli Abruzzi 24,}\\{\footnotesize 10129 Torino, Italy. Member of INdAM research group GNCS.}
	}
	
	\maketitle

	\begin{abstract}
		A numerical method for coupled 3D-1D problems with discontinuous solutions at the interfaces is derived and discussed. This extends a previous work on the subject where only continuous solutions were considered. Thanks to properly defined function spaces a well posed 3D-1D problem is obtained from the original fully 3D problem and the solution is then found by a PDE-constrained optimization reformulation. This is a domain decomposition strategy in which unknown interface variables are introduced and a suitably defined cost functional, expressing the error in fulfilling interface conditions, is minimized constrained by the constitutive equations on the subdomains. The resulting discrete problem is robust with respect to geometrical complexity thanks to the use of independent discretizations on the various subdomains. Meshes of different sizes can be used without affecting the conditioning of the discrete linear system, and this is a peculiar aspect of the considered formulation. An efficient resolution strategy is further proposed, based on the use of a gradient based solver and yielding a method ready for parallel implementation. A numerical experiment on a problem with known analytical solution shows the accuracy of the method, and two examples on more complex configurations are proposed to address the applicability of the approach to practical problems.
	\end{abstract}
\subsection*{Keywords}
3D-1D coupling; optimization based domain-decomposition; non conforming mesh\\
AMS Subject Classification: 65N30, 65N50, 68U20, 86-08

\section{Introduction}
The present work deals with simulations in domains with embedded cylindrical, or nearly cylindrical inclusions, with cross-section sizes much smaller than their length and than the domain scale. This kind of problem is typical of a large variety of applications, ranging from the study of living tissues, where the inclusions are constituted, e.g., by the capillaries or by the vessels of the lymphatic system  \cite{Cattaneo2014,Notaro2016,Koppl2020}, to the study of the interaction between plant roots and the soil \cite{Schroder2012,Koch2018}, or fibre-reinforced materials \cite{Steinbrecher2020,Llau2016}, and geological applications \cite{Gjerde2018a,Gjerde2020, CerLauZun2019}. In these cases, it is often convenient to treat the inclusions as one dimensional objects, thus actually collapsing the cross-sections on the centrelines, avoiding the complexity and the overhead of generating good quality meshes in the interior of the small inclusions. Such geometrical reduction is however non-trivial from a mathematical standpoint as it ends up in a 3D problem which is coupled to problems on the 1D domains, and non standard function spaces need to be used to allow the definition of a well posed trace operator for 3D functions on 1D manifolds. 

An analysis of 3D problems with singular source terms is available in Refs.~\cite{Dangelo2012,DangeloQuarteroni2008}, where the solution is placed in suitable weighted Sobolev spaces. In Ref.~\cite{Dangelo2012} a finite element based approach is used for the resolution of the problem and optimal convergence trends for the error are observed in the used spaces. Problems with singular source terms are also considered in Ref.~\cite{Tornberg2004}, where Dirac delta sources are replaced by regularizing terms, compactly supported. Regularizing techniques are suggested also in Ref.~\cite{Heltai2019}, and in Ref.~\cite{Koch2020} where line source terms are approximated by suitable kernel functions that distribute the source in a three dimensional neighbourhood of the line. A splitting technique is proposed in Ref.~\cite{Gjerde2019}, where the solution is seen as the sum of a known low regularity term and a regular correction term that is computed solving an elliptic problem with source term and boundary data depending on a chosen extension operator for the singular source. In Ref.~\cite{Koppl2018} a lifting technique is employed for the singular source term, whereas a domain decomposition approach based on the use of Lagrange multipliers is suggested in Ref.~\cite{Kuchta2021}. The derivation of a coupled 3D-1D problem, starting from an original equidimensional formulation is proposed in Ref.~\cite{Zunino2019}, under suitable assumptions on the behavior of the solution in the inclusions that allow the definition of averaging operators.

Here an extension of the approach described in Ref.~\cite{BGS3D1D2022} for 3D-1D coupled elliptic problems is proposed, allowing to deal with discontinuous solution at the interfaces. This is relevant for a large variety of practical applications, as e.g. the description of biological tissues, where the boundary between the 3D and the 1D domains is a semi-permeable membrane and the Starling equation applies \cite{Cattaneo2014,Zunino2019,Koch2020,Kuchta2021}. A well posed 3D-1D coupled problem is derived from the original 3D-3D problem through an appropriate choice of the functional space for the solution, which allows to define extension and trace operators between spaces on 3D and 1D domains. Further, a domain decomposition approach is employed, introducing additional interface unknowns to decouple the 3D problem from the problems on the 1D inclusions, and a cost functional is designed such that the solution is obtained as the minimum of the functional constrained by the constitutive equations on the subdomains. The original work in Ref.~\cite{BGS3D1D2022} deals with continuous solution at the boundary between the 3D and the 1D domains. Here, while keeping the same structure of this original method, different interface variables are introduced for the domain decomposition process, which result in a novel setting for the PDE constrained optimization problem described. The problem is discretized using finite elements on non conforming meshes, thus providing a great flexibility in the choice of the meshes that can be independently defined on each subdomain, and a numerical scheme suitable for parallel computing is obtained thanks to the use of domain decomposition combined to a gradient scheme for the resolution of the resulting discrete problem. 

The problem considered is presented in Section~\ref{not_and_form} in equi-dimensional form, and in Section~\ref{var_form} its 3D-1D formulation is derived in weak form in ad-hoc function spaces. The problem is re-written as a PDE constrained optimization problem in Section~\ref{PDE_constr}, and presented in discrete form in Section~\ref{Discrete}. The approach used to solve the obtained linear system is discussed in Section~\ref{solve_strategy}, whereas some numerical tests are reported in Section~\ref{Num_res} and finally conclusions are reported in Section~\ref{Concl}.

\section{Notation and formulation of the fully 3D coupled problem}\label{not_and_form}

Let $\Omega {\subset} \mathbb{R}^3$ be a convex domain in which a generalized cylinder $\Sigma {\subset} \mathbb{R}^3$ is embedded. We denote by $\Lambda=\left\lbrace\bm{\lambda}(s), s \in (0,S) \right\rbrace$ the centerline of $\Sigma$, while $\Gamma=\left\lbrace \Gamma(s), ~s \in [0,S]\right\rbrace $ is the lateral surface of $\Sigma$. In the following we assume, for simplicity, that $\Lambda$ is a rectilinear segment in the three-dimensional space. We denote by $\Sigma(s)$ the transversal sections of the cylinder as $s$ ranges in the interval $[0,S]$ and  by $\Gamma(s)$ their boundary. We suppose the sections to have an elliptic shape, with $R(s)$ being the maximum axes length of the ellipse centered in $\bm{\lambda}(s)$. For the two extreme sections of the cylinder we adopt the compact notation  $\Sigma_0=\Sigma(0)$ and $\Sigma_S=\Sigma(S)$. For the derivation of the model problem we assume that $\Sigma_0$ and $\Sigma_S$ lie on the boundary $\partial \Omega$, but the extension to more general cases is straightforward. The portion of $\Omega$ that does not include the cylinder is denoted by $D=\Omega \setminus \Sigma$. We define $\partial D^e=\partial \Omega \setminus \left\lbrace \Sigma_0 \cup \Sigma_S\right\rbrace$, referring to it as the \textit{external boundary} of $D$, with $\partial D =\partial D^e \cup \Gamma$. In case the extreme sections of $\Sigma$ were inside $\Omega$, $\partial D^e$ would coincide with $\partial \Omega$.

Let us now consider the following diffusion problem, with unknown pressures $u$ in $D$ and $\tilde{u}$ in $\Sigma$:

\begin{center}
	
	\begin{minipage}{0.47\textwidth}
		\begin{center}
			\textbf{3D-problem on} $\bm{D}$:
			\begin{align}
			&-\nabla \cdot (K \nabla u)=f & \text{in } D\label{eqOmega}\\
			&u=0 &\text{on } \partial D^e \label{cbOmega}\\
			&K\nabla u \cdot \bm{n}=\beta(\tilde u_{|_{\Gamma}} -u_{|_{\Gamma}}) \hspace{-0.3cm}&\text{ on } \Gamma \label{phi_u}
			\end{align}		
		\end{center}
	\end{minipage}
	\begin{minipage}{0.49\textwidth}
		\begin{center}
			\textbf{3D-problem on} $\bm{\Sigma}$:
			\begin{align}
			&-\nabla \cdot({\tilde{K}}\nabla\tilde{u})=g&~\text{ in } \Sigma \label{eqSigma}\\
			&\tilde{u}=0 &\text{ on } \Sigma_0 \cup \Sigma_S \label{cbSigma}\\
			&\tilde{K}\nabla \tilde{u} \cdot \bm{\tilde{n}}=\beta(u_{|_{\Gamma}}-\tilde u_{|_{\Gamma}})  \hspace{-0.9cm}&\text{ on } \Gamma\label{phi_utilde}	
			\end{align}	
		\end{center}
	\end{minipage}
	
\end{center}
\medskip

\noindent Vectors $\bm{n}$ and $\bm{\tilde{n}}=-\bm{n}$ are the outward pointing unit normal vectors to $\Gamma$, respectively for $D$ and $\Sigma$; $K$, $\tilde{K}$ and $\beta$ are positive scalars, while $f$ and $g$ denote source terms. For the sake of simplicity we consider homogeneous Dirichlet boundary conditions on $\Sigma_0$ and $\Sigma_S$ and on $\partial D ^e$. 
Equations \eqref{phi_u} and \eqref{phi_utilde} allow us to couple the two problems imposing flux conservation. According to these equations, the flux  across $\Gamma$ is directly proportional to the jump of the pressures, with $\beta$ denoting the permeability coefficient of the membrane $\Gamma$. Different coupling conditions could be considered, for example adding a pressure continuity constraint and consequently not linking the flux definition to the pressure jump, as done in Ref.~\cite{BGS3D1D2022}. The choice of the interface condition depends of course on the properties of the interface, and thus on the kind of application.

Let us now suppose that $R$ is much smaller than the size of $\Omega$ and than the longitudinal length $L$ of the cylinder itself, in particular. This allows us to assume that the variables defined on $\Sigma$ or on $\Gamma$ are actually only functions of the coordinate $s$, considering negligible their variation on the cross-sections of the inclusion. 
This is the key point that allows us, in the next section, to work out a well-posed 3D-1D coupled problem from equations \eqref{eqOmega}-\eqref{phi_utilde}.

\section{Variational formulation of the 3D-1D problem} \label{var_form}
A 3D-1D coupled problem is obtained from problem \eqref{eqOmega}-\eqref{phi_utilde}, after writing it in variational form in suitable function spaces, as here described. Let us, thus, introduce the spaces  
\begin{align*}
&H_0^1(D)=\left\lbrace v \in H^1(D) : v_{|_{\partial D^e}}=0 \right\rbrace,\\ 
&H_0^1(\Sigma)=\left\lbrace v \in H^1(\Sigma): v_{|_{\Sigma_0}}=v_{|_{\Sigma_S}}=0\right\rbrace, \\
&H_0^1(\Lambda)=\left\lbrace v \in H^1(\Lambda): v(0)=v(S)=0\right\rbrace,
\end{align*}
a trace operator
\begin{equation*}
\gamma_{_\Gamma}:H^1(D)\cup H^1(\Sigma)\rightarrow H^{\frac{1}{2}}(\Gamma),\text{ such that }\gamma_{_\Gamma}v=v_{|_\Gamma} ~\forall v \in H^1(D)\cup H^1(\Sigma)
\end{equation*} 
and two extension operators 
\begin{equation*}
\mathcal{E}_{_\Sigma}: H^1(\Lambda) \rightarrow H^1(\Sigma) ~\text{  and  }~\mathcal{E}_{_\Gamma}: H^1(\Lambda) \rightarrow H^{\frac{1}{2}}(\Gamma)
\end{equation*}
such that, given $\hat{v} \in H_0^1(\Lambda)$, $\mathcal{E}_{_\Sigma}\hat{v}(s)$ and $\mathcal{E}_{_\Gamma}\hat{v}(s)$ are the uniform extension of $\hat{v}(s)$ respectively to $\Sigma(s)$ and to $\Gamma(s)$. We observe that $\mathcal{E}_{_\Gamma}=\gamma_{_\Gamma}\circ \mathcal{E}_{_\Sigma}$.
Once denoted by $\hat{V}$ the space $H_0^1(\Lambda)$, let us further consider the spaces:
\begin{align*}
&\widetilde{V}=\lbrace v \in H_0^1(\Sigma): v =\exts\hat{v}, ~\hat{v} \in \hat{V} \rbrace, \\
&\mathcal{H}^{\Gamma}=\lbrace v \in H^{\frac{1}{2}}(\Gamma): v =\extg\hat{v}, ~\hat{v} \in \hat{V} \rbrace,\\
&V_D=\left\lbrace v \in H_0^1(D): \gamma_{_\Gamma}v \in \mathcal{H}^{\Gamma}\right\rbrace,
\end{align*} 
such that functions in $\widetilde{V}$ and in $\hsp$ are respectively the uniform extension to $\Sigma$ and $\Gamma$ of functions in $\hat{V}$ and functions in $V_D$ have trace on $\Gamma$ belonging to $\hsp$. Functions in such spaces fit the assumptions we have made on the negligible variation of the variables on the cross-sections of $\Sigma$ and $\Gamma$. Denoting by $(\cdot,\cdot)_{X}$ the scalar product on a generic space $X$, the variational formulation of problem \eqref{eqOmega}-\eqref{phi_utilde} can be written as:
\textit{find} $(u, \tilde{u}) \in V_D \times \widetilde{V}$ \textit{such that}
\begin{align}
&({K}\nabla u, \nabla v)_{L^2(D)}-\left( \beta(\trg\tilde u -\trg u),\gamma_{_\Gamma}v \right)_{\mathcal{H}^{\Gamma}}=(f,v)_{L^2(D)}~&\forall v \in V_D \label{coup_eq1}\\
&({\tilde{K}}\nabla \tilde{u}, \nabla \tilde{v})_{L^2(\Sigma)}+\left( \beta(\trg \tilde u -\trg u),\gamma_{_\Gamma}v \right)_{\mathcal{H}^{\Gamma}}=(g,\tilde{v})_{{L^2(\Sigma)}} &\forall \tilde{v} \in \widetilde{V} \label{coup_eq2}
\end{align}
Let us introduce two auxiliary variables $\psid,\psis \in \hsp$, in order to formally decouple the two equations. Denoting by $X'$ the dual of the generic space $X$, the problem is thus rewritten as: \textit{find} $(u, \tilde{u}) \in V_D \times \widetilde{V}$, $\psid \in \hsp$ and $\psis \in \hsp$ \textit{such that}
\begin{align}
&({K}\nabla u, \nabla v)_{L^2(D)}+\left( \beta \trg u,\gamma_{_\Gamma}v \right) _{\mathcal{H}^{\Gamma}}-\left( \beta\psis ,\gamma_{_\Gamma}v \right) _{\mathcal{H}^{\Gamma}}=(f,v)_{{L^2(D)}}~&\forall v \in V_D\label{eq_var_u} \\
&({\tilde{K}}\nabla \tilde{u}, \nabla \tilde{v})_{L^2(\Sigma)}+\left(  \beta\trg \tilde u,\gamma_{_\Gamma}\tilde{v} \right)_{\mathcal{H}^{\Gamma}}-\left(  \beta\psid,\gamma_{_\Gamma}\tilde{v} \right)_{\mathcal{H}^{\Gamma}}=(g,\tilde{v})_{L^2(\Sigma)} &\forall \tilde{v} \in \widetilde{V} \label{eq_var_uhat}
\end{align}\vspace{-0.7cm}
\begin{align}
&\left\langle \gamma_{_\Gamma}u-\psid,\eta\right\rangle_{\mathcal{H}^{\Gamma},{\mathcal{H}^{\Gamma}}'}=0&~\forall \eta \in  {\mathcal{H}^{\Gamma}}' \label{condpsi_u}\\
&\left\langle \gamma_{_\Gamma}\tilde{u}-\psis,\eta\right\rangle_{\mathcal{H}^{\Gamma},{\mathcal{H}^{\Gamma}}'}=0&~\forall \eta \in  {\mathcal{H}^{\Gamma}}'\label{condpsi_hat}.
\end{align}
Let us remark that Equations \eqref{eq_var_u}-\eqref{eq_var_uhat} could also be written as
\begin{align*}
&({K}\nabla u, \nabla v)_{L^2(D)}+\left( \beta\psid,\gamma_{_\Gamma}v \right) _{\mathcal{H}^{\Gamma}}-\left( \beta\psis ,\gamma_{_\Gamma}v \right) _{\mathcal{H}^{\Gamma}}=(f,v)_{L^2(D)}~&\forall v \in V_D \\
&({\tilde{K}}\nabla \tilde{u}, \nabla \tilde{v})_{L^2(\Sigma)}+\left(  \beta\psis,\gamma_{_\Gamma}\tilde{v} \right)_{\mathcal{H}^{\Gamma}}-\left(  \beta\psid,\gamma_{_\Gamma}\tilde{v} \right)_{\mathcal{H}^{\Gamma}}=(g,\tilde{v})_{L^2(\Sigma)} &\forall \tilde{v} \in \widetilde{V}. 
\end{align*}
However, formulation \eqref{eq_var_u}-\eqref{eq_var_uhat} is preferred, as it allows to have an empty Dirichlet boundary on either $\partial D^e$ or $\Sigma_0,\Sigma_s$. This is a desired property for domain decomposition purposes.

Thanks to the adopted functional spaces, problem \eqref{eq_var_u}-\eqref{condpsi_hat} can be easily reduced to a 3D-1D coupled problem. Let us observe that, given $\eta \in {\hsp}'$ and $\rho \in \hsp$
\begin{equation*}
\left\langle \rho,\eta \right\rangle_{{\mathcal{H}^{\Gamma}}, {\mathcal{H}^{\Gamma}}'}=\int_{\Gamma}\rho\eta~d\Gamma=\int_0^S\Big( \int_{\Gamma(s)}\rho\eta~dl\Big) ds.
\end{equation*}
Since $\rho \in \hsp$, there exists $\hat{\rho} \in \hat{V}$ such that $\extg\hat{\rho}=\rho$ and thus $\int_{\Gamma(s)}\rho  ~dl=|\Gamma(s)|\hat{\rho}(s)$. Hence we can introduce $\overline{\eta}\in \hat{V}'$ such that 
\begin{equation*}
\int_0^S\Big( \int_{\Gamma(s)}\rho\eta~dl\Big) ds=\int_0^S|\Gamma(s)|\hat{\rho}(s)\overline{\eta}(s)~ds=\left\langle \hat{\rho},|\Gamma|\overline{\eta}\right\rangle_{\hat{V}',\hat{V}},
\end{equation*}
where $|\Gamma(s)|$ is the section perimeter size at $s\in[0,S]$. Similarly, if we consider $\rho,w \in \hsp$, then
\begin{equation*}
\left(  \rho,w \right)_{\mathcal{H}^{\Gamma}}=\int_0^S|\Gamma(s)|\hat{\rho}(s)\hat{w}(s)~ds=\left(|\Gamma|\hat{\rho},\hat{w}\right)_{L^2(\Lambda)}
\end{equation*}
with $\extg\hat{\rho}=\rho$ and $\extg\hat{w}=w$ . Finally let us observe that
\begin{equation*}
({\tilde{K}}\nabla \tilde{u}, \nabla \tilde{v})_{L^2(\Sigma)}=\int_{\Sigma}{\tilde{K}}\nabla \tilde{u}\nabla\tilde{v}~d\sigma=\int_0^S{\tilde{K}}|\Sigma(s)|\cfrac{d\hat{u}}{ds}~\cfrac{d\hat{v}}{ds}~ds
\end{equation*}
where $\hat{u},\hat{v} \in \hat{V}$ are such that $\tilde{u}=\exts\hat{u}$, $\tilde{v}=\exts\hat{v}$ and $|\Sigma(s)|$ is the section area at $s \in [0,S]$.
{Problem  (\ref{eq_var_u})-(\ref{condpsi_hat}) can now be rewritten as a reduced 3D-1D coupled problem:}
\textit{Find $(u,\hat{u}) \in {V_D}\times\hat{V}$, $\hpsid\in \hat{V}$ and $\hpsis \in \hat{V}$ such that:}
\begin{align}
&({{K}}\nabla u, \nabla v)_{L^2(D)}+\left( |\Gamma|\beta\check{u},\check{v}\right) _{L^2(\Lambda)}-\left( |\Gamma|\beta\hpsis,\check{v} \right) _{L^2(\Lambda)}=(f,v)_{L^2(D)} \label{equaz1} \\[-0.5em]  &\hspace{7cm}\forall v \in {V_D}, \check{v} \in \hat{V}: \gamma_{_\Gamma} v=\mathcal{E}_{_\Gamma}\check{v} \nonumber \\
&\Big( {\tilde{K}}|\Sigma|\frac{d\hat{u}}{ds},\frac{d\hat{v}}{ds}\Big)_{{L^2(\Lambda)}}\hspace{-0.2cm}+\left( |\Gamma|\beta\hat{u},\hat{v} \right) _{L^2(\Lambda)}-\left(  |\Gamma|\beta\hpsid,\hat{v}\right) _{L^2(\Lambda)}=(|\Sigma|\overline{\overline{g}},\hat{v})_{{L^2(\Lambda)}}\label{equaz2} \\[-0.5em]  &\hspace{10cm} \forall \hat{v} \in \hat{V}\nonumber
\end{align}\vspace{-0.5cm}
\begin{align}
&\left\langle |\Gamma|(\check{u}-\hpsid),\overline{\eta}\right\rangle_{\hat{V},\hat{V}'}=0  &\gamma_{_\Gamma}u=\extg\check{u},\forall \overline{\eta} \in \hat{V}'\label{condiz1}\\
&\left\langle |\Gamma|(\hat{u}-\hpsis),\overline{\eta}\right\rangle_{\hat{V},\hat{V}'}=0 &\forall \overline{\eta} \in \hat{V}'\label{condiz2}
\end{align}
with $\overline{\overline{g}}(s)=\frac{1}{|\Sigma(s)|}\int_{\Sigma(s)}g~d\sigma$, being $g$ sufficiently regular.

\section{PDE-constrained optimization problem}\label{PDE_constr}
Conditions (\ref{condiz1}) and (\ref{condiz2}) can be replaced by the minimization of a cost functional mimicking the error committed in the fulfillment of such constraints. At this aim let us define
\begin{eqnarray}
J(\hpsid,\hpsis)&=&\cfrac{1}{2}\left( ||\gamma_{_\Gamma}u(\hpsis)-\psid||_{\hsp}^2+||\gamma_{_\Gamma}\tilde{u}(\hpsid)-\psis||_{\hsp}^2\right) \nonumber \\
&=&\cfrac{1}{2}\left( ||\gamma_{_\Gamma}u(\hpsis)-\extg \hpsid||_{\hsp}^2+||\gamma_{_\Gamma}\exts \hat{u}(\psid)-\extg \hpsis||_{\hsp}^2\right)	\label{functional}
\end{eqnarray}
to be minimized constrained by \eqref{equaz1} and \eqref{equaz2}. In order to work out the PDE-constrained optimization formulation of the problem in a compact form, let us define the linear operators
$A: {V_D} \rightarrow {V_D'}$, $\widehat{A}:\hat{V} \rightarrow \hat{V}'$, $S:\hat{V}\rightarrow V_D'$ and $\widehat{D}:\hat{V}\rightarrow \hat{V}'$ such that, for any $u \in {V_D}, ~\check{u} \in \hat{V}: \trg u=\extg\check{u}$, $\hat{u} \in \hat{V}$, $\hpsis,\hpsid \in \hat{V}$:
\begin{align}
&\left\langle Au,v\right\rangle_{{V_D'},{V_D}}=({K}\nabla u, \nabla v)_{{L^2(D)}}+\left( |\Gamma|\beta\check{u},\check{v}\right)_{L^2(\Lambda)}\label{A} \\[-0.5em] &\hspace{7cm} v \in {V_D}, ~\check{v} \in \hat{V}: \trg v=\extg\check{v}\nonumber
\end{align}\vspace{-0.7cm}
\begin{align}
&\left\langle \widehat{A}\hat{u},\hat{v}\right\rangle_{\hat{V}',\hat{V}}=\Big({\tilde{K}}|\Sigma|\cfrac{d\hat{u}}{ds},\cfrac{d\hat{v}}{ds}\Big)_{L^2(\Lambda)}+\left(  |\Gamma|\beta\hat{u},\hat{v}\right)_{L^2(\Lambda)}& \hat{v} \in \hat{V} \label{Ahat}\\
&\left\langle S \hpsis,v \right\rangle _{V_D',V_D}=(|\Gamma|\beta\hpsis,\check{v})_{L^2(\Lambda)} & \hspace{-2cm}v \in {V_D}, ~\check{v} \in \hat{V} :\trg v=\extg\check{v} \label{C}\\
&\left\langle \widehat{D} \hpsid,\hat{v} \right\rangle _{\hat{V}',\hat{V}}=(|\Gamma| \beta\hpsid,\hat{v})_{L^2(\Lambda)}& \hat{v} \in \hat{V}. \label{Chat}
\end{align}
The respective adjoints will be denoted as $A^*: {V_D} \rightarrow {V_D'}$, $\widehat{A}^*:\hat{V} \rightarrow \hat{V}'$, $S:V_D \rightarrow \hat{V}'$ and $\widehat{D}^*:\hat{V} \rightarrow \hat{V}'$. If we further define 
\begin{align}
&F\in {V_D'} \text{ s.t. } F(v)=(f,v)_{{L^2(D)}}, &v \in {V_D}\\
&G\in\hat{V}' \text{ s.t. }G(\hat{v})=(|\Sigma|\overline{\overline{g}},\hat{v})_{L^2(\Lambda)}, &\hat{v} \in \hat{V},\label{G}
\end{align}
equations \eqref{equaz1}-\eqref{equaz2} can be written as: 
\begin{align}
&Au-S\hpsis=F \label{eq1}\\
&\widehat{A}\hat{u}-\widehat{D}\hpsid=G.\label{eq2}
\end{align}
Finally, the PDE-constrained optimization problem can be written as
\begin{equation}
\min_{\hpsid,\hpsis\in \hat{V}}J(\hpsid,\hpsis) \text{ subject to } \eqref{eq1}-\eqref{eq2} \label{minJ}
\end{equation}

We now provide some results on the optimal control and the stepsize of the steepest descent method for Problem \eqref{minJ}.

\begin{prop}
	Let us consider the trace operator $\gamma_{_\Gamma}:{V_D}\rightarrow \hsp$ and the extension operators $\mathcal{E}_{_\Sigma}:\hat{V} \rightarrow \tilde{V}$ and $\mathcal{E}_{_\Gamma}=\gamma_{_\Gamma}\circ \mathcal{E}_{_\Sigma}:\hat{V}\rightarrow \hsp$, whose respective adjoints are $\gamma_{_\Gamma}^*:{\hsp}'\rightarrow {V_D'}$, ${\exts}^*:\tilde{V}'\rightarrow \hat{V}'$ and $\extg^*:{\hsp}'\rightarrow \hat{V}'$ and let $\Theta_{\hat{V}}:\hat{V}\rightarrow\hat{V}'$ and $\Theta_{\hsp}:\hsp\rightarrow {\hsp}'$ be  Riesz isomorphisms. Then the optimal control $(\hpsid,\hpsis)$ that provides the solution to \eqref{minJ} is such that
	\begin{align}
	& \Theta_{\hat{V}}^{-1}(\extg^*\Theta_{\hsp}(\extg\hpsid-\gamma_{_\Gamma}u(\hpsis))+\widehat{D}^*\hat{p})=0\\
	&\Theta_{\hat{V}}^{-1}(\extg^*\Theta_{\hsp}(\extg\hpsis-\extg\hat{u}(\hpsid))+S^*p)=0
	\end{align}
	where $p \in {V_D}$ and $\hat{p} \in \hat{V}$ are the solutions respectively to
	\begin{align}
	&A^*p=\gamma_{_\Gamma}^*\Theta_{\hsp}(\gamma_{_\Gamma}u(\hpsis)-\extg\hpsid)\label{p}\\
	&\widehat{A}^* \hat{p}=\extg^*\Theta_{\hsp}(\extg\hat{u}(\hpsid)-\extg\hpsis)\label{phat}
	\end{align}
	\label{prop1}
\end{prop}
\begin{proof}
	Let us compute the Fréchet derivatives of $J$ with respect to the control variables $\hpsid$ and $\hpsis$. To this end, we introduce the increments $\delta\hpsid, \delta \hpsis \in \hat{V}$ and we recall that, for $\star=D,\Sigma$, there exists $ \delta {\psi}_\star \in \mathcal{H}^{\Gamma}~:~\delta \psi_\star=\mathcal{E}_{_\Gamma}\delta\hat{ \psi}_\star$. We have:
	\begin{align*}
	&\cfrac{\partial J}{\partial \hpsid}\big(\hpsid+\delta\hpsid,\hpsis\big) =\left(\gamma_{_\Gamma}{u}(\hpsis)-\psid,-\delta\psid\right)_{\hsp}+\left(\gamma_{_\Gamma}\tilde{u}(\hpsid)-\psis,\gamma_{_\Gamma}\tilde{u}(\delta\hpsid)\right)_{\hsp}\\
	&=\left( \mathcal{E}_{_\Gamma}\hpsid-\gamma_{_\Gamma}u(\hpsis),\extg\delta \hpsid\right)_{\hsp}+\left(\trg \exts \hat{u}(\hpsid)-\extg \hpsis,\trg \exts \hat{u}(\delta\hpsid) \right)_{\hsp}=\\
	&=\left\langle \extg^*\Theta_{\hsp}(\extg\hpsid-\gamma_{_\Gamma}u(\hpsis)),\delta \hpsid\right\rangle_{\hat{V}',\hat{V}}+\left\langle\extg^*\Theta_{\hsp}(\extg\hat{u}(\hpsid)-\extg\hpsis),\hat{u}(\delta\hpsid)\right\rangle_{\hat{V}',\hat{V}}\\
	&=\left\langle \extg^*\Theta_{\hsp}(\extg\hpsid-\gamma_{_\Gamma}u(\hpsis)),\delta \hpsid\right\rangle_{\hat{V}',\hat{V}}+\left\langle \widehat{A}^*\hat{p},\widehat{A}^{-1}\widehat{D}\delta \hpsid \right\rangle_{\hat{V}',\hat{V}}=\\
	&=\left\langle \extg^*\Theta_{\hsp}(\extg\hpsid-\gamma_{_\Gamma}u(\hpsis)),\delta \hpsid\right\rangle_{\hat{V}',\hat{V}} +\left\langle \widehat{D}^*\hat{p},\delta \hpsid \right\rangle_{{\hat{V}}',\hat{V}}=\\
	&=\left( \Theta_{\hat{V}}^{-1}(\extg^*\Theta_{\hsp}(\extg\hpsid-\gamma_{_\Gamma}u(\hpsis))+\widehat{D}^*\hat{p}),\delta \hpsid\right)_{{L^2(\Lambda)}};
	\end{align*}
	\begin{align*}
	&\cfrac{\partial J}{\partial \hpsis}(\hpsid,\hpsis+\delta \hpsis)=\left(\gamma_{_\Gamma}{u}(\hpsis)-\psid,\gamma_{_\Gamma}u(\delta\hpsis)\right)_{\hsp}+\left(\gamma_{_\Gamma}\tilde{u}(\hpsid)-\psis,-\delta\psis\right)_{\hsp}\\
	&=\left\langle  \trg^*\Theta_{\hsp}(\gamma_{_\Gamma}u(\hpsis)-\extg\hpsid),u(\delta\hpsis)\right\rangle _{V_D',V_D}+\left( \extg\hpsis-\gamma_{_\Gamma}\exts\hat{u}(\hpsid),\extg\delta\hpsis\right)_{\hsp}=\\[0.5em]
	&=\left\langle A^*p,A^{-1}S\delta \hpsis \right\rangle_{{V_D'},{V_D}}+\left\langle \extg^*\Theta_{\hsp}(\extg\hpsis-\extg\hat{u}(\hpsid)),\delta\hpsis\right\rangle_{\hat{V}',\hat{V}}=\\[0.5em]
	&=\left\langle S^*p,\delta \hpsis \right\rangle_{{\hat{V}}',\hat{V}} +\left\langle \extg^*\Theta_{\hsp}(\extg\hpsis-\extg\hat{u}(\hpsid)),\delta\hpsis\right\rangle_{\hat{V}',\hat{V}}=\\
	&=\left( \Theta_{\hat{V}}^{-1}(S^*p+\extg^*\Theta_{\hsp}(\extg\hpsis-\extg\hat{u}(\hpsid)),\delta \hpsis\right)_{L^2(\Lambda)},
	\end{align*}
	which yield the thesis.
\end{proof}

\noindent From the derivatives computed in Proposition \ref{prop1}, we now define the quantities
\begin{align}
&\delta\hpsid=\Theta_{\hat{V}}^{-1}(\extg^*\Theta_{\hsp}(\extg\hpsid-\gamma_{_\Gamma}u(\hpsis))+\widehat{D}^*\hat{p})\in \hat{V}\\
&\delta \hpsis=\Theta_{\hat{V}}^{-1}(\extg^*\Theta_{\hsp}(\extg\hpsis-\extg\hat{u}(\hpsid))+S^*p) \in \hat{V}.
\end{align}
Then the following proposition holds:
\begin{prop}
	Given the variable $\mathcal{X}$, let us increment it by a step $\zeta \delta \mathcal{X}$, where $\delta \mathcal{X}=(\delta \hpsid, \delta \hpsis)$. The steepest descent method corresponds to the stepsize 
	$\zeta=-\frac{\mathcal{N}}{\mathcal{D}}$ with
	\begin{align*}
	&\mathcal{N}=\left( \delta \hpsid,\delta \hpsid\right)_{L^2(\Lambda)}+\left( \delta \hpsis,\delta \hpsis\right)_{L^2(\Lambda)}\\
	&\mathcal{D}=\left\langle S\delta\hpsis,\delta p \right\rangle_{{V_D'},{V_D}}-\left\langle\extg^*\Theta_{\hsp}(\trg \delta u-\extg \delta\hpsid),\delta \hpsid\right\rangle_{\hat{V}',\hat{V}}+\left\langle \widehat{D}\delta\hpsid,\delta \hat{p} \right\rangle_{\hat{V}',\hat{V}}+\\[-0.5em]&\qquad -\left\langle\extg^*\Theta_{\hsp}(\extg \delta\hat{u}-\extg \delta\hpsis),\delta \hpsis\right\rangle_{\hat{V}',\hat{V}}
	\end{align*}
	and where
	\begin{align*}
	&\delta u=u(\delta \hpsis)=A^{-1}S\delta \hpsis \in {V_D},\\
	&\delta \hat{u}=\hat{u}(\delta \hpsid )=\widehat{A}^{-1}\widehat{D}\delta\hpsid \in \hat{V}
	\end{align*}
	and $\delta p \in {V_D}$, $\delta \hat{p} \in \hat{V}$ are such that:
	\begin{align*}
	&A^*\delta p=\gamma_{_\Gamma}^*\Theta_{\hsp}(\gamma_{_\Gamma}\delta u-\extg\delta\hpsid)\\
	&\widehat{A}^* \delta\hat{p}=\extg^*\Theta_{\hsp}(\extg\delta \hat{u}-\extg\delta \hpsis)
	\end{align*}	
	\label{prop2}
\end{prop}

\begin{proof}
	It is sufficient to set to zero the derivative $\frac{\partial J(\mathcal{X}+\zeta \delta \mathcal{X})}{\partial \zeta}$. In the computation that follows we adopt the lighter notation:
	\begin{equation*}
	u=u(\hpsis);~~\delta u=u(\delta\hpsis);\qquad
	\hat{u}=\hat{u}(\hpsid);~~\delta \hat{u}=\hat{u}(\delta\hpsid).
	\end{equation*}
	\begin{align*}
	&J(\mathcal{X}+\zeta\delta \mathcal{X})=J(\hpsid+\zeta\delta \hpsid,\hpsis+\zeta\delta\hpsis)=\\[0.5em]
	&=\cfrac{1}{2}\left(\gamma_{_\Gamma}u(\hpsis+\zeta\delta\hpsis)-\psid-\zeta\delta\psid,\gamma_{_\Gamma}u(\hpsis+\zeta\delta\hpsis)-\psid-\zeta\delta \psid \right)_{\hsp}\\[-0.2em]
	&+\cfrac{1}{2}\left(\gamma_{_\Gamma}\tilde{u}(\hpsid+\zeta\delta \hpsid)-\psis-\zeta\delta\psis,\gamma_{_\Gamma}\tilde{u}(\hpsid+\zeta\delta \hpsid)-\psis-\zeta\delta \psis \right)_{\hsp}=\\[0.5em]
	&=\cfrac{1}{2}\left(\gamma_{_\Gamma} u+\zeta\gamma_{_\Gamma}\delta u-\extg\hpsid-\zeta\extg\delta\hpsid,\gamma_{_\Gamma}u+\zeta\gamma_{_\Gamma}\delta u-\extg\hpsid-\zeta\extg\delta\hpsid\right)_{\hsp}+\\[-0.2em]
	&+\cfrac{1}{2}\left(\gamma_{_\Gamma}\exts\hat{u}+\zeta\gamma_{_\Gamma}\exts\delta\hat{u}-\extg\hpsis-\zeta\extg\delta\hpsis,\gamma_{_\Gamma}\exts\hat{u}+\zeta\gamma_{_\Gamma}\exts\delta\hat{u}-\extg\hpsis-\zeta\extg\delta\hpsis\right)_{\hsp}\\[0.5em]	
	&=J(\hpsid,\hpsis)+\zeta\left( \gamma_{_\Gamma}u-\extg\hpsid,\gamma_{_\Gamma}\delta u-\extg\delta\hpsid\right)_{\hsp}+\zeta\left( \extg(\hat{u}-\hpsis),\extg(\delta\hat{u}-\delta\hpsis)\right)_{\hsp}\\[-0.2em]
	&+ \cfrac{\zeta^2}{2}\left( \gamma_{_\Gamma}\delta u-\extg\delta\hpsid,\gamma_{_\Gamma}\delta u-\extg\delta\hpsid\right) _{\hsp}+\cfrac{\zeta^2}{2}\left( \extg(\delta\hat{u}-\delta\hpsis),\extg(\delta \hat{u}-\delta\hpsis)\right) _{\hsp}
	\end{align*}
	\begin{align*}
	&\cfrac{\partial J(\mathcal{X}+\zeta \delta \mathcal{X})}{\partial \zeta}=\left( \gamma_{_\Gamma}u-\extg\hpsid,\gamma_{_\Gamma}\delta u-\extg\delta\hpsid\right)_{\hsp}+\left( \extg(\hat{u}-\hpsis),\extg(\delta\hat{u}-\delta\hpsis)\right)_{\hsp}+\\[-0.2em]
	&\quad+\zeta\left( \gamma_{_\Gamma}\delta u-\extg\delta\hpsid,\gamma_{_\Gamma}\delta u-\extg\delta\hpsid\right) _{\hsp}+\zeta\left( \extg(\delta\hat{u}-\delta\hpsis),\extg(\delta\hat{u}-\delta\hpsis)\right) _{\hsp}=0
	\end{align*}
	\begin{align*}
	\Rightarrow ~\zeta=-\cfrac{\left( \gamma_{_\Gamma}u-\extg\hpsid,\gamma_{_\Gamma}\delta u-\extg\delta\hpsid\right)_{\hsp}+\left( \extg\hat{u}-\extg\hpsis,\extg\delta\hat{u}-\extg\delta\hpsis\right)_{\hsp}}{\left( \gamma_{_\Gamma}\delta u-\extg\delta\hpsid,\gamma_{_\Gamma}\delta u-\extg\delta\hpsid\right) _{\hsp}+\left( \extg\delta\hat{u}-\extg\delta\hpsis,\extg\delta\hat{u}-\extg\delta \hpsis\right) _{\hsp}}
	\end{align*}
	Rearranging properly the terms we get $\zeta=-\frac{\mathcal{N}}{\mathcal{D}}$ with 
	\begin{align*}
	\mathcal{N}&=\left\langle A^*p,A^{-1}S\delta\hpsis\right\rangle _{{V_D'},{V_D}}-\left\langle\extg^*\Theta_{\hsp}(\trg u-\extg \hpsid),\delta \hpsid\right\rangle_{\hat{V}',\hat{V}}+\\[-0.5em]
	&\qquad +\left\langle\widehat{A}^*\hat{p},\widehat{A}^{-1}\widehat{D}\delta\hpsid \right\rangle_{\hat{V}',\hat{V}}-\left\langle\extg^*\Theta_{\hsp}(\extg \hat{u}-\extg \hpsis),\delta \hpsis\right\rangle_{\hat{V}',\hat{V}}=\\
	&=\left\langle S^*p,\delta \hpsis \right\rangle_{{\hat{V}}',\hat{V}}-\left\langle\extg^*\Theta_{\hsp}(\trg u-\extg \hpsid),\delta \hpsid\right\rangle_{\hat{V}',\hat{V}}+\left\langle \widehat{D}^*\hat{p},\delta \hpsid \right\rangle_{{\hat{V}}',\hat{V}}+\\[-0.5em] &\qquad -\left\langle\extg^*\Theta_{\hsp}(\extg \hat{u}-\extg \hpsis),\delta \hpsis\right\rangle_{\hat{V}',\hat{V}}=\\	
	&=\left( \delta \hpsid,\delta \hpsid\right)_{L^2(\Lambda)}+\left( \delta \hpsis,\delta \hpsis\right)_{L^2(\Lambda)}
	\end{align*}
	and
	\begin{align*}
	\mathcal{D}&=\left\langle A^{^{-1}}S\delta\hpsis,A^{^*}\delta p\right\rangle_{{V_D},{V_D'}}-\left\langle\extg^*\Theta_{\hsp}(\trg \delta u-\extg \delta\hpsid),\delta \hpsid\right\rangle_{\hat{V}',\hat{V}}+\\[-0.5em]&\quad+\left\langle \widehat{A}^{^{-1}}\widehat{D}\delta\hpsid,\widehat{A}^*\delta \hat{p}\right\rangle_{\hat{V},\hat{V}'} -\left\langle\extg^*\Theta_{\hsp}(\extg \delta\hat{u}-\extg \delta\hpsis),\delta \hpsis\right\rangle_{\hat{V}',\hat{V}}=\\
	&=\left\langle S\delta\hpsis,\delta p \right\rangle_{{V_D'},{V_D}}-\left\langle\extg^*\Theta_{\hsp}(\trg \delta u-\extg \delta\hpsid),\delta \hpsid\right\rangle_{\hat{V}',\hat{V}}+\left\langle \widehat{D}\delta\hpsid,\delta \hat{p} \right\rangle_{\hat{V}',\hat{V}}+\\[-0.5em]
	&\quad -\left\langle\extg^*\Theta_{\hsp}(\extg \delta\hat{u}-\extg \delta\hpsis),\delta \hpsis\right\rangle_{\hat{V}',\hat{V}}
	\end{align*}
	that yields the thesis.
\end{proof}

\section{Discrete matrix formulation}\label{Discrete}
\newcommand{\I}{\mathcal{I}}
In this section we work out the discrete matrix formulation of problem \eqref{minJ}.
In general, the 3D-1D coupling does not present particular issues in the discrete framework. Nonetheless our approach has the additional advantage of allowing for the use of non conforming meshes: thanks to the optimization framework, the partitions of the 1D inclusions can be defined in a completely independent manner from the surrounding 3D mesh, without any theoretical or practical constraint on mesh sizes. Further, the proposed formulation provides the direct computation of interface variables, and it allows to decouple the 3D problem from the 1D problems, thus paving the way to the use of possibly different constitutive equations and to the application of efficient solvers based on parallel computing techniques.

For the sake of generality, we consider $\mathcal{I}$ segments of different length and orientation crossing the domain $\Omega$. The segments are defined as $\Lambda_i=\left\lbrace \bm{\lambda}_i(s),s \in (0, S_i)\right\rbrace $, $i=1,...,\I$ and they represent the centerlines of $\I$ cylindrical inclusions $\Sigma_i$. 

The proposed approach can easily handle intersections among inclusions centrelines. 
Intersecting segments are split into sub-segments in correspondence of their intersection point $q$. In this way, $q$ always corresponds to a segment endpoint, in which pressure continuity and flux conservation are constrained. It is to remark that a variety of intersection modes is possible for the original three dimensional inclusions. As an example, 3D inclusions might partially overlap whereas the corresponding centrelines might not intersect. By considering here only intersections between centrelines, we implicitly assume that the intersection volume of the corresponding three dimensional inclusions is small and can be reduced to a point in the scale of the domain. A deeper investigation on the treatment of different intersection models is out of the scope of the present work.   

After having extended the domain D to the whole $\Omega$, let us consider a tetrahedral mesh $\mathcal{T}$ of domain $\Omega$, on which we define Lagrangian finite element basis functions $\left\lbrace \varphi_k \right\rbrace_{k=1}^{N}$, such that $U=\sum_{k=1}^{N}U_k\varphi_k$ is the discrete approximation of pressure $u$. On each segment $\Lambda_i$ we build three different partitions, independent from each other and from $\mathcal{T}$. We denote them by $\hat{\mathcal{T}_i}$, $\tau^{D}_i$ and $\tau^{\Sigma}_i$ and we define the basis functions $\left\lbrace\hat{\varphi}_{i,k} \right\rbrace _{k=1}^{\hat{N}_i}$ on $\hat{\mathcal{T}_i}$, $\left\lbrace \theta_{i,k}^D\right\rbrace_{k=1}^{N_i^{D}}$ on $\tau^{D}_i$ and $\left\lbrace \theta_{i,k}^\Sigma\right\rbrace_{k=1}^{N_i^{\Sigma}}$ on $\tau^{\Sigma}_i$, with $\hat{N}_i$, $N_i^{D}$ and $N_i^{\Sigma}$ denoting the number of DOFs of the discrete approximations of the variables $\hat{u}_i$, $\hat{\psi}_{_D,i}$ and $\hat{\psi}_{_\Sigma,i}$ respectively. Such approximations are defined as
\begin{equation*}
\hat{U}_i=\sum_{k=1}^{\hat{N}_i}\hat{U}_{i,k}~\hat{\varphi}_{i,k}, \quad \Psi^D_i=\sum_{k=1}^{N_i^{D}}\Psi_{i,k}^D~\theta_{i,k}^D, \quad \Psi^\Sigma_i=\sum_{k=1}^{N_i^{\Sigma}}\Psi_{i,k}^\Sigma ~\theta_{i,k}^\Sigma.
\end{equation*}
We then define the following matrices:
\begin{align*}
&\bm{A} \in \mathbb{R}^{N\times N} \text{ s.t. } (A)_{kl}=\int_{\Omega}K\nabla\varphi_k\nabla\varphi_l ~d\omega+\sum_{i=1}^{\I}\int_{\Lambda_i}\beta_i|\Gamma_i(s)|{\varphi_k}_{|_{\Lambda_i}}{\varphi_l}_{|_{\Lambda_i}} ds\\[0.8em]
&\bm{\hat{A}_i} \in \mathbb{R}^{\hat{N}_i\times \hat{N}_i} \text{ s.t. } (\hat{A}_i)_{kl}=\int_{\Lambda_i}\tilde{K}_i|\Sigma_i(s)|\frac{d\hat{\varphi}_{i,k}}{ds}\frac{d\hat{\varphi}_{i,l}}{ds} ~ds+\int_{\Lambda_i}\beta_i|\Gamma_i(s)|\hat{\varphi}_{i,k}\hat{\varphi}_{i,l}~ds
\end{align*}
\begin{align*}
&\bm{\hat{D}_i}^{\beta} \in \mathbb{R}^{\hat{N}_i\times N_i^{D}} \text{ s.t. } (\hat{D}_i^\beta)_{kl}=\int_{\Lambda_i}\beta_i|\Gamma_i(s)|{\hat{\varphi}_{i,k}~\theta_{i,l}^D}~ds\\
&\bm{S_i}^{\beta} \in \mathbb{R}^{N\times N_i^{\Sigma}} \text{ s.t. } (S_i^{\beta})_{kl}=\int_{\Lambda_i}\beta_i|\Gamma_i(s)|{\varphi_k}_{|_{\Lambda_i}}\theta_{i,l}^{\Sigma}~ds,
\end{align*}
and the vectors 
\begin{equation*}
f\in \mathbb{R}^N \text{ s.t. } f_k=\int_{\Omega}f\varphi_k~d\omega,\qquad \qquad
g_i\in \mathbb{R}^{\hat{N}_i} \text{ s.t. } (g_i)_k=\int_{\Lambda_i}|\Sigma_i(s)|\overline{\overline{g}}~\hat{\varphi}_{i,k}~ds.
\end{equation*}
Setting $\hat{N}=\sum_{i=1}^{\I}\hat{N}_i$, $N_D=\sum_{i=1}^{\I}N_i^{D}$ and $N_\Sigma=\sum_{i=1}^{\I}\ N_i^{\Sigma}$, 
we can group the matrices as follows:
\begin{align*}
\bm{\hat{D}}^{\beta}=\text{diag}\left( \bm{\hat{D}_1}^{\beta},...,\bm{\hat{D}_{\I}^{\beta}}\right)   \in \mathbb{R}^{\hat{N}\times N_D} \qquad
\bm{S}^{\beta}=\left[ \bm{S_1}^{\beta}, \bm{S_2}^{\beta},...,\bm{S_{\I}}^{\beta}\right] \in \mathbb{R}^{N\times N_\Sigma}
\end{align*}
and
\begin{displaymath}
\bm{\hat{A}}=\left[
\begin{array}{cc}
\text{diag}\left(  \bm{\hat{A}_1},...,\bm{\hat{A}_{\I}}\right)  & \bm{Q}^T\\
\bm{Q} & \bm{0}
\end{array}\right]=\left[\begin{array}{cc}
\bm{\hat{A}^{\sharp}} & \bm{Q}^T\\
\bm{Q} & \bm{0}
\end{array}
\right]
\end{displaymath}
where matrix $\bm{Q}$ simply equates the DOFs at the extrema of connected sub-segments and allow us to enforce continuity through Lagrange multipliers. Let us observe how $\bm{\hat{A}}=\bm{\hat{A}^\sharp}$ in case no intersections occur among segments.
Finally we can write
\begin{align}
&\bm{A}U-\bm{S}^\beta\Psi_{\Sigma}=f\label{eq1discr}\\
&\bm{\hat{A}}\hat{U}-\bm{\hat{D}}^\beta\Psi_{D}=g\label{eq2discr}
\end{align}
with
\begin{align*}
&\hat{U}=\left[\hat{U}_1^T,...,\hat{U}_{\I}^T \right]^T \in \mathbb{R}^{\hat{N}}; \quad g=[g_1^T,g_2^T,...,g_{\I}^T]^T\in\mathbb{R}^{\hat{N}}\\ 
&\Psi_D=\left[(\Psi_1^D)^T,...,(\Psi_\I^D)^T \right]^T \in \mathbb{R}^{N_D}; \quad \Psi_\Sigma=\left[(\Psi_1^\Sigma)^T,...,(\Psi_\I^\Sigma)^T \right]^T\in \mathbb{R}^{N_\Sigma}.
\end{align*}
The discrete functional is derived from \eqref{functional} replacing the norms in $\hsp$ with norms in $L^2(\Lambda)$ and summing over the $\mathcal{I}$ inclusions. First we define matrices
\begin{align*}
&\bm{G_i} \in \mathbb{R}^{N \times N} \text{ s.t. } (G_i)_{kl}=\int_{\Lambda_i}{\varphi_k}_{|_{\Lambda_i}}{\varphi_l}_{|_{\Lambda_i}}ds,\\
&\bm{\hat{G}_i} \in \mathbb{R}^{\hat{N}_i \times \hat{N}_i} \text{ s.t. } (\hat{G}_i)_{kl}=\int_{\Lambda_i}{\hat{\varphi}}_{i,k}~{\hat{\varphi}}_{i,l}~ds,\\
&\bm{{M}_i^{D}} \in \mathbb{R}^{N_i^{D} \times N_i^{D}}  \text{ s.t. } ({M}_i^{D})_{kl}=\int_{\Lambda_i}\theta_{i,k}^D~\theta_{i,l}^D~ds,\\
&\bm{{M}_i^{\Sigma}} \in \mathbb{R}^{N_i^{\Sigma} \times N_i^{\Sigma}}  \text{ s.t. } ({M}_i^{\Sigma})_{kl}=\int_{\Lambda_i}\theta_{i,k}^\Sigma~\theta_{i,l}^\Sigma~ds,\\
&\bm{D_i} \in \mathbb{R}^{N\times N_i^{D}} \text{ s.t. } (D_i)_{kl}=\int_{\Lambda_i}{{\varphi_{k}}_{|_{\Lambda_i}}\theta_{i,l}^D}~ds,\\ &\bm{\hat{S}_i} \in \mathbb{R}^{\hat{N}_i\times N_i^{\Sigma}} \text{ s.t. } (\hat{S_i})_{kl}=\int_{\Lambda_i}{\hat{\varphi}}_{i,k}~\theta_{i,l}^\Sigma~ds,
\end{align*}
and then
\begin{equation*}
\label{Gdef}
\bm{G}=\sum_{i=1}^{\I}\bm{G}_i \in \mathbb{R}^{N \times N} \qquad \bm{\hat{G}}=\text{diag}\left( \bm{\hat{G}_1},...,\bm{\hat{G}_{\I}}\right)  \in \mathbb{R}^{\hat{N} \times\hat{N}}
%\end{bmatrix}
\end{equation*}
\begin{equation*}
\bm{M^D}=\text{diag}\left(\bm{M_1^{D}},...,\bm{M_{\I}^{D}} \right) \in \mathbb{R}^{N_{D}\times N_{D}}, \qquad  \bm{M^{\Sigma}}=\text{diag}\left(\bm{M_1^{\Sigma}},...,\bm{M_{\I}^{\Sigma}} \right) \in \mathbb{R}^{N_{\Sigma}\times N_{\Sigma}}
\end{equation*}
\begin{align*}
&\bm{D}=\left[ \bm{D_1}, \bm{D_2},...,\bm{D_{\I}}\right]  \in \mathbb{R}^{N\times N_D} \qquad  \bm{\hat{S}}=\text{diag}\left(  \bm{\hat{S}_1},...,\bm{\hat{S}_{\I}}\right)   \in \mathbb{R}^{\hat{N}\times N_\Sigma}
\end{align*}
The discrete cost functional then reads:
\begin{align}
\tilde{J}=\cfrac{1}{2}\Big( U^T&\bm{G}U-U^T\bm{D}\Psi_D-\Psi_D^T\bm{D}^TU+\Psi_D^T\bm{M^D}\Psi_D+\nonumber\\&+\hat{U}^T\bm{\hat{G}}\hat{U}-\hat{U}^T\bm{\hat{S}}\Psi_\Sigma -\Psi_\Sigma^T\bm{\hat{S}}^T\hat{U}+\Psi_\Sigma^T\bm{M^\Sigma}\Psi_\Sigma\Big)\label{Jtilde} 
\end{align}
Finally, the discrete matrix formulation of the 3D-1D problem can be written as:
\begin{align}
\min_{(\Psi_D,\Psi_\Sigma)}\tilde{J}(\Psi_D,\Psi_\Sigma) \text{ subject to } \eqref{eq1discr}-\eqref{eq2discr} \label{minJtilde}
\end{align}
First order optimality conditions for problem \eqref{minJtilde} are collected in the saddle-point system
\begin{equation}
\bm{\mathcal{K}}\hspace{-0.12cm}=\hspace{-0.12cm}
\begin{bmatrix}
\bm{G} & \bm{0} & -\bm{D} & \bm{0} &\bm{A}^T &\bm{0} \\
\bm{0} & \bm{\hat{G}} & \bm{0} &-\bm{\hat{S}} &\bm{0} &\bm{\hat{A}}^T \\ 
-\bm{D}^T & \bm{0} & \bm{M^D} & \bm{0} & \bm{0} & (-\bm{{\hat{D}}^{\beta}})^T \\
\bm{0} &-\bm{\hat{S}}^T & \bm{0} & \bm{M^\Sigma} & (-\bm{S^\beta})^T &\bm{0}\\
\bm{A} & \bm{0} & \bm{0} & -\bm{S^\beta} &\bm{0} &\bm{0}\\
\bm{0} &\bm{A}^T & -\bm{{\hat{D}}^{\beta}} & \bm{0} & \bm{0} &\bm{0}
\end{bmatrix}\hspace{-0.12cm}; \quad  
\bm{\mathcal{K}}
\begin{bmatrix}
U\\\hat{U}\\ \Psi_D \\ \Psi_\Sigma \\ -P\\ -\hat{P}
\end{bmatrix}\hspace{-0.15cm}=\hspace{-0.15cm}\begin{bmatrix}
0 \\0  \\ 0 \\ 0 \\f \\g
\end{bmatrix}\label{KKT}
\end{equation}

\begin{prop}
	\label{discretewellposedness}
	Matrix $\bm{\mathcal{K}}$ in \eqref{KKT} is non-singular and the unique solution of \eqref{KKT} is equivalent to the solution of the optimization problem  \eqref{minJtilde}.
\end{prop}
The following lemma is used to prove Proposition~\ref{discretewellposedness}.

\begin{lemma}\label{lemma1}
	Let matrix $\bm{\mathcal{A}}\in \mathbb{R}^{(N+\hat{N})\times(N+\hat{N}+N_D+N_\Sigma)}$ be as
	\begin{displaymath}
	\bm{\mathcal{A}}=
	\begin{bmatrix}
	\bm{A} & \bm{0} & \bm{0} & -\bm{S^\beta}\\
	\bm{0} &\bm{A}^T & -\bm{{\hat{D}}^{\beta}} & \bm{0}
	\end{bmatrix}
	\end{displaymath}
	and let $\bm{\mathcal{G}}\in \mathbb{R}^{(N+\hat{N}+N_D+N_\Sigma)\times(N+\hat{N}+N_D+N_\Sigma)}$ be defined as
	\begin{displaymath}
	\bm{\mathcal{G}}=
	\begin{bmatrix}
	\bm{G} & \bm{0} & -\bm{D} & \bm{0}\\
	\bm{0} & \bm{\hat{G}} & \bm{0} &-\bm{\hat{S}} \\ 
	-\bm{D}^T & \bm{0} & \bm{M^D} & \bm{0}  \\
	\bm{0} &-\bm{\hat{S}}^T & \bm{0} & \bm{M^\Sigma}
	\end{bmatrix}.
	\end{displaymath}
	Then matrix $\bm{\mathcal{A}}$ is full rank and matrix $\bm{\mathcal{G}}$ is symmetric positive definite on $\ker(\bm{\mathcal{A}})$.
\end{lemma}

\begin{proof}
	%%%%%%%%%%%%%%%%%%%%%%%%%%%%%%%%%%%%%%%%%%
	Matrix $\bm{\mathcal{A}}$ is full rank for the ellipticity of operators $A$ in \eqref{A} and $\hat{A}$ in \eqref{Ahat}, whereas matrix $\bm{\mathcal{G}}$ is symmetric positive semi-definite as
	\begin{displaymath}
	\bm{\Phi}=\begin{bmatrix} U\\\hat{U}\\ \Psi_D \\ \Psi_\Sigma  \end{bmatrix}, \quad
	\bm{\Phi}^T \bm{\mathcal{G}} 
	\bm{\Phi} = \frac12\sum_{i=1}^\I \left(\| U_{|\Lambda_i} - \Psi_i^D \|^2_{L^2(\Lambda_i)} + \| \hat{U}_i - \Psi_i^\Sigma \|^2_{L^2(\Lambda_i)}\right)\geq 0
	\end{displaymath} 
	Let $e_k=\begin{bmatrix} e_k^D \\ e_k^\Sigma \end{bmatrix} $ be the $k$-th element of the canonical basis of $\mathbb{R}^{N_D+N_{\Sigma}}$ and let $z_k\in\ker{(\bm{\mathcal{A}})}$ be defined as:
	\begin{displaymath}
	z_k=\begin{bmatrix}
	\bm{A}^{-1}(-\bm{S^\beta} e_k^\Sigma)\\
	\bm{\hat{A}}^{-1}(-\bm{\hat{D}^{\beta}} e_k^D)\\
	e_k
	\end{bmatrix} := \begin{bmatrix} z_k^D \\ z_k^\Sigma \\ e_k^D \\ e_k^\Sigma \end{bmatrix}
	\end{displaymath}
	Thus it is either $e_k^D \neq \bm{0}$ and $z_k^\Sigma \neq \bm{0}$, either $e_k^\Sigma \neq \bm{0}$ and $z_k^D \neq \bm{0}$, and consequently 
	\begin{displaymath}
	z_k^T \bm{\mathcal{G}} z_k = \frac12 \left(\| z_k^D - e_k^D \|^2_{L^2(\Lambda)} + \| z_k^\Sigma - e_k^\Sigma \|^2_{L^2(\Lambda)}\right) > 0,
	\end{displaymath}
	for a certain segment $\Lambda$, depending on $k$.
	As a consequence, $z_k\not\in\ker{(\bm{\mathcal{G}})}$ for any $k=1,\ldots, N^D+N^\Sigma$. The vector space $\ker{(\bm{\mathcal{A}})}=\text{span}\{z_1,\ldots,z^{N^D+N^\Sigma}\}$ is a subspace of $\text{Im}(\bm{\mathcal{G}})$, and $\ker{(\bm{\mathcal{G}})} \cap \ker{(\bm{\mathcal{A}})} =\{\bm{0} \}$.
\end{proof}

The proof of Proposition~\ref{discretewellposedness} derives from classical arguments of quadratic programming, observing that 
\begin{displaymath}
\bm{\mathcal{K}}=\begin{bmatrix}
\bm{\mathcal{G}} & \bm{\mathcal{A}}^T\\
\bm{\mathcal{A}} & \bm{0}
\end{bmatrix}.
\end{displaymath}

\section{Solving strategies}\label{solve_strategy}
Solving system \eqref{KKT} is equivalent to solve the optimum problem \eqref{minJtilde}. A different resolution approach is however proposed, based on an iterative solver and allowing to take full advantage of the decoupling introduced by the proposed method.

Let us formally replace in the cost functional \eqref{Jtilde} the expressions $U=\bm{A}^{-1}(\bm{S^\beta}\Psi_{\Sigma}+f)$ and $\hat{U}=\bm{\hat{A}}^{-1}(\bm{\hat{D}^\beta}\Psi_{D}+g)$ and let us set $\mathcal{X}=[\Psi_D^T,\Psi_\Sigma^T]^T$, obtaining
\begin{align}
&J^\star(\Psi_D,\Psi_\Sigma)=\cfrac{1}{2}\Big( (\bm{A}^{-1}\bm{S^\beta}\Psi_{\Sigma}+\bm{A}^{-1}f)^T\bm{G}(\bm{A}^{-1}\bm{S^\beta}\Psi_{\Sigma}+\bm{A}^{-1}f)+\nonumber\\
&\qquad-(\bm{A}^{-1}\bm{S^\beta}\Psi_{\Sigma}+\bm{A}^{-1}f)^T\bm{D}\Psi_D-\Psi_D^T\bm{D}^T(\bm{A}^{-1}\bm{S^\beta}\Psi_{\Sigma}+\bm{A}^{-1}f)+\nonumber\\
&\qquad+\Psi_D^T\bm{M^D}\Psi_D+(\bm{\hat{A}}^{-1}\bm{\hat{D}^\beta}\Psi_D+\bm{\hat{A}}^{-1}g)^T\bm{\hat{G}}(\bm{\hat{A}}^{-1}\bm{\hat{D}^\beta}\Psi_D+\bm{\hat{A}}^{-1}g)+\nonumber\\
&\qquad
-(\bm{\hat{A}}^{-1}\bm{\hat{D}^\beta}\Psi_D+\bm{\hat{A}}^{-1}g)^T\bm{\hat{S}}\Psi_\Sigma-\Psi_\Sigma^T\bm{\hat{S}}^T(\bm{\hat{A}}^{-1}\bm{\hat{D}^\beta}\Psi_D+\bm{\hat{A}}^{-1}g)+\nonumber\\
&\qquad+\Psi_\Sigma^T\bm{M^\Sigma}\Psi_\Sigma\Big)=\nonumber\\
&=\cfrac{1}{2}~\mathcal{X}^T\begin{bmatrix}
(\bm{\hat{D}^\beta})^T\bm{\hat{A}}^{-T}\bm{\hat{G}}\bm{\hat{A}}^{-1}\bm{\hat{D}^\beta}+\bm{M^D} &\quad -\bm{D}^T\bm{A}^{-1}\bm{S^\beta}-(\bm{\hat{D}^\beta})^T\bm{\hat{A}}^{-T}\bm{\hat{S}}\\\\
-(\bm{S^\beta})^T\bm{A}^{-T}\bm{D}-\bm{\hat{S}}^T\bm{\hat{A}}^{-1}\bm{\hat{D}^\beta} &(\bm{S^\beta})^T\bm{A}^{-T}\bm{G}\bm{A}^{-1}\bm{S^\beta}+\bm{M^\Sigma}
\end{bmatrix}\mathcal{X}\nonumber\\[0.5em]
&\qquad+ \begin{bmatrix}g^T
\bm{\hat{A}}^{-T}\bm{\hat{G}}\bm{\hat{A}}^{-1}\bm{\hat{D}^\beta}-f^T\bm{A}^{-T}\bm{D} \\\\ f^T\bm{A}^{-T}\bm{G}\bm{A}^{-1}\bm{S^\beta}-g^T\bm{\hat{A}}^{-T}\bm{\hat{S}}
\end{bmatrix}\mathcal{X}+\nonumber\\
&\qquad+\cfrac{1}{2}\left( f^T\bm{A}^{-T}\bm{G}\bm{A}^{-1}f+g^T\bm{\hat{A}}^{-T}\bm{\hat{G}}\bm{\hat{A}}^{-1}g\right)= \nonumber \\ 
&=\cfrac{1}{2}\left( \mathcal{X}^T\bm{\mathcal{M}}\mathcal{X}+2d\mathcal{X}+q\right).\label{Jcompact}
\end{align}
Matrix $\bm{\mathcal{M}}$ is symmetric positive definite, given the equivalence of this formulation with the saddle-point system \eqref{KKT}. This allows us to perform the minimization of the unconstrained functional \eqref{Jcompact} via a gradient based scheme, looking for the minimum as the solution of
\begin{equation}
\nabla J^*=\bm{\mathcal{M}}\mathcal{X}+d=0. \label{systCg}
\end{equation}
A preconditioner $\bm{\mathcal{P}}$ can be defined for the resolution of system \eqref{systCg}. In particular we set 
\begin{equation}
\bm{\mathcal{P}}=\begin{bmatrix}
(\bm{\hat{D}^\beta})^T(\bm{\hat{A}^\sharp})^{-T}\bm{\hat{G}}(\bm{\hat{A}}^\sharp)^{-1}\bm{\hat{D}^\beta}+\bm{M^D} &\bm{0}\\\\
\bm{0} &\bm{M^\Sigma}
\end{bmatrix}. \label{prec}
\end{equation}
where $\bm{\hat{A}^\sharp}=\text{diag}\left(  \bm{\hat{A}_1},...,\bm{\hat{A}_{\I}}\right)$. This means that the top-left block of matrix $\bm{\mathcal{P}}$  corresponds exactly to the top-left block of matrix $\bm{\mathcal{M}}$ in case no intersections among the segment occur. Otherwise, $(\bm{\hat{A}^\sharp})^{-1}$ is an approximation of the inverse of matrix  $\bm{\hat{A}}$ which can be built inverting independently the matrices related to the single segments and which maintains a block-diagonal structure, i.e. $(\bm{\hat{A}}^\sharp)^{-1}=\text{diag}\left(  \bm{\hat{A}_1}^{-1},...,\bm{\hat{A}_{\I}}^{-1}\right)$. For what concerns the bottom-right block, only matrix $\bm{M^{\Sigma}}$ is kept with respect to the same block of matrix $\bm{\mathcal{M}}$, so that even this portion of the preconditioner can be built block-diagonalwise, assembling matrices which are independently related to each single 1D inclusion.

The conjugate gradient scheme which is employed to solve system $\eqref{systCg}$ is reported in Algorithm \ref{grad_congAlgo}. The quantity $\bm{\mathcal{M}}\delta \mathcal{X}$, whose computation is required at each iteration in the algorithm, can actually be obtained without explicitly building matrix $\bm{\mathcal{M}}$. In fact, after some computations we obtain
\begin{equation*}
\bm{\mathcal{M}}\delta\mathcal{X}=\begin{bmatrix}
(\bm{\hat{D}^\beta})^T\delta \hat{P}-\bm{D}^T\delta U+\bm{M}^D\delta \Psi_D\\
(\bm{S}^\beta)^T\delta {P}-\bm{\hat{S}}^T\delta \hat U+\bm{M}^\Sigma\delta \Psi_\Sigma
\end{bmatrix}
\end{equation*}
where $\delta \mathcal{X}=[\delta \Psi_D^T,\delta\Psi_\Sigma^T]$, and $\delta U$, $\delta \hat{U}$, $\delta P$, $\delta \hat{P}$ are the solutions of the linear systems
\begin{align*}
& \bm{A}\delta U=\bm{S}^\beta\delta \Psi_\Sigma \qquad 
& \bm{\hat{A}}\delta \hat{U}=\bm{\hat{D}}^\beta\delta \Psi_D\\
&\bm{A}\delta P=\bm{G}\delta U-\bm{D}\delta \Psi_D &\bm{\hat{A}}\delta\hat{ P}=\bm{\hat{G}}\delta \hat{U}-\bm{\hat{S}}\delta \Psi_\Sigma
\end{align*}
which require the resolution of local sub-problems on the 1D segments and on the 3D domain.

\begin{algorithm}
	\SetAlgoLined
	Guess $\mathcal{X}_0=[\Psi_{D,0}^T, \Psi_{\Sigma,0}^T]^T$ \\
	$r_0=\bm{\mathcal{M}}\mathcal{X}_0+d$;\\
	Solve $\bm{\mathcal{P}}z_0=r_0$;\\
	set $\delta \mathcal{X}_0=-z_0$ and $k=0$;\\
	\While{$\frac{||r_k||}{||d||} >toll$ \label{toll}}{
		$\zeta_k=\cfrac{r_k^Tz_k}{\delta \mathcal{X}_{k}^T\bm{\mathcal{M}}\delta \mathcal{X}_{k}}$;\\
		$\mathcal{X}_{k+1}=\mathcal{X}_k+\zeta_k\delta \mathcal{X}_{k}$;\\
		$r_{k+1}=r_k+\zeta_k\bm{\mathcal{M}}\delta \mathcal{X}_{k}$;\\
		Solve $\bm{\mathcal{P}}z_{k+1}=r_{k+1}$;\\
		$\beta_{k+1}=\cfrac{r_{k+1}^Tz_{k+1}}{r_{k}^Tz_{k}}$;\\
		$\delta \mathcal{X}_{k+1}=-z_{k+1}+\beta_{k+1}\delta \mathcal{X}_k$;\\
		$k=k+1$;}
	\caption{Conjugate gradient method for $\bm{\mathcal{M}}\mathcal{X}+d=0$}
	\label{grad_congAlgo}
\end{algorithm}

\section{Numerical results}\label{Num_res}
In this section we present three numerical examples to better highlight the characteristics of the proposed approach. The simulations are performed using linear finite elements on the 3D and 1D non-conforming meshes, independently generated on the sub-domains. Parameter $h$ denotes the maximum diameter of the tetrahedra for the 3D mesh, while other three parameters, namely  $\hat{\delta}_{u,i}$, $\delta_{D,i}$ and $\delta_{\Sigma,i}$ express the refinement level of the 1D meshes $\hat{\mathcal{T}_i}$, $\tau^{D}_i$, $\tau^{\Sigma}_i$, $i=1,...\I$, respectively. Each of these three parameters represents the ratio between the number of nodes in the 1D mesh and the number of intersections of the segment $\Lambda_i$ with the faces of the tetrahedra in $\mathcal{T}$. 
In the simulations, for simplicity, we adopt unique, but possibly different, values of  $\hat{\delta}_{u}$, $\delta_{D}$ and $\delta_{\Sigma}$ for the various segments: for this reason we drop, in the following, the segment index $i=1,...\I$ for these parameters. In all cases, linear Lagrangian finite elements on tethrahedra are used in the 3D domains and piecewise continuous linear basis functions on equally spaced meshes are chosen for the 1D functions.

\subsection[Test1]{Test Problem 1 (TP1)}
\begin{figure}
	\centering
	\includegraphics[width=0.55\linewidth]{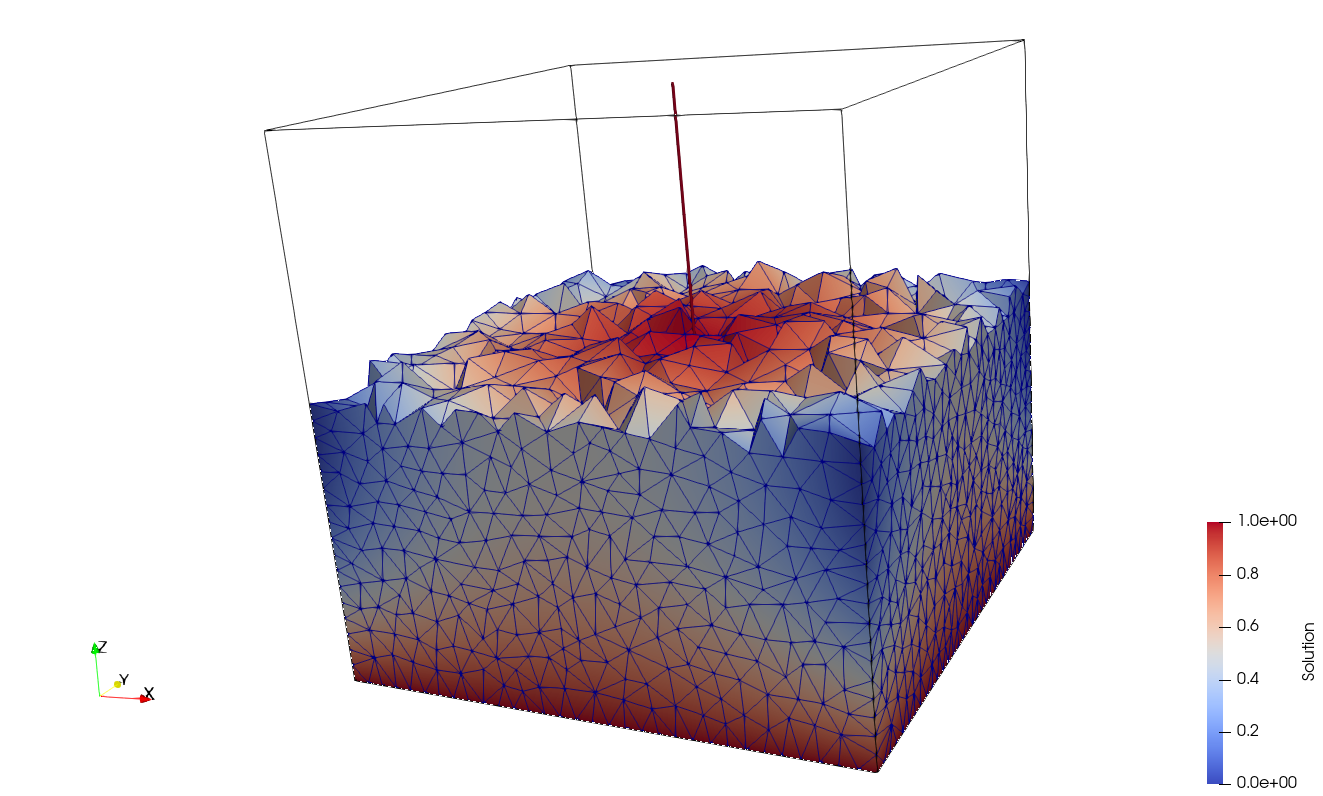}
	\includegraphics[width=0.40\linewidth]{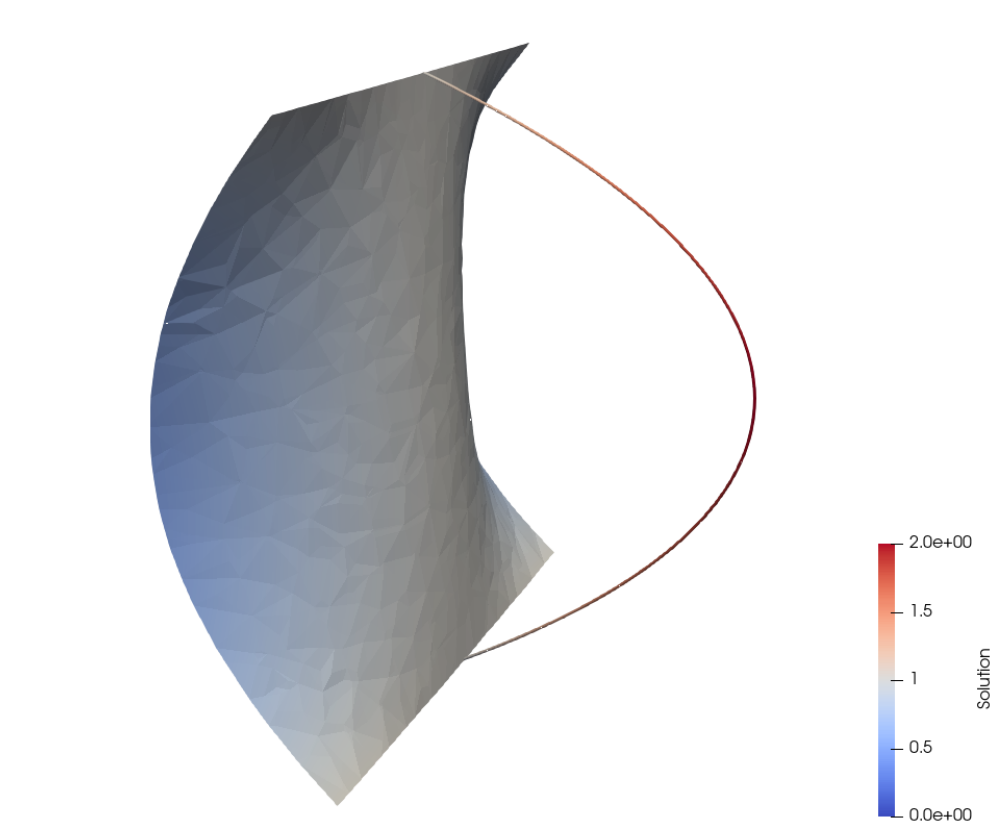}
	\caption{TP1: left: view of the numerical solution inside the cube; right: solution obtained on the segment and on a section of the cube parallel to the $z$-axis and containing $\Lambda$. Parameters $h = 0.13$, $\hat{\delta}_u = 1$, $\delta_D=\delta_{\Sigma} = 0.5$.}
	\label{section_view}
\end{figure}

Let us consider a cube $\Omega$ of edge $l=2$ centered in the axes origin and whose faces are parallel to the coordinate axes. Let us further consider a cylinder $\Sigma$ of radius $\hat{R}=10^{-2}$ and height $h=2$ whose centreline $\Lambda$ lies on the $z$ axis (see Figure~\ref{section_view}, on the left). Let us denote by $\partial \Omega_{l}$, $\partial \Omega_{+}$ and $\partial \Omega_{-}$ respectively the lateral, the top and the bottom faces of the cube.
We aim at solving a problem in the form of \eqref{equaz1}-\eqref{condiz2}, obtained by reducing $\Sigma$ to its centerline, with the following data:
\begin{align*}
&K=1, \qquad f(x,y,z)=2-x^2-y^2-2z^2,\\ &\tilde{K}(z)=\frac{z^2}{3}+\frac{1}{2},\qquad \overline{\overline{g}}=3\\
&\beta=\frac{2\hat{R}}{2+\hat{R}^2}
\end{align*}
The problem is completed with appropriate boundary conditions such that the exact solution is:
\begin{align}
&u_{ex}(x,y,z)=\frac{1}{2}(x^2+y^2)(z^2-1)+1& \text{ in } \Omega \label{uestest1}\\
&\hat{u}_{ex}(z)=2-z^2&\label{uhatestest1} \text{ on } \Lambda.
\end{align} 
In particular we consider Neumann boundary conditions on $\partial\Omega_{+}$ and $\partial \Omega_{-}$, whereas Dirichlet boundary conditions are imposed on $\partial\Omega_l$. Dirichlet boundary conditions equal to 1 are imposed at segment endpoints.

The obtained solution is shown in Figure~\ref{section_view}: on the left, the 3D solution is shown on a portion of the domain, whereas, on the right, the solution $U$ on the $y-z$ plane containing the $z$-axis is plotted along with the solution $\hat{U}$, with solution values reported along the $x$-axis. The computational mesh used for this solution has parameters $h=0.083$, $\hat{\delta}_u=1$, $\delta_{D}=0.5$ and $\delta_{\Sigma}=0.5$, corresponding to $N=4155$ DOFs in the cube and $\hat{N}=57$ DOFs on the segment.

\begin{figure}
	\centering
	\includegraphics[width=0.46\linewidth]{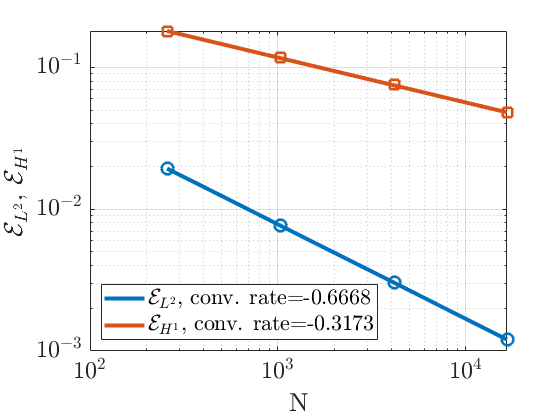}%
	\includegraphics[width=0.46\linewidth]{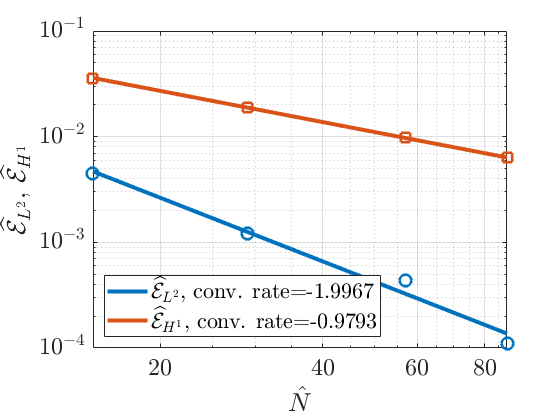}
	\caption{TP1: trend of the $L^2$ and $H^1$-norms of the relative errors under mesh refinement. On the left: error on the cube with respect to (\ref{uestest1}); on the right: error on the segment with respect to \eqref{uhatestest1}. Other parameters: $\hat{\delta}_{u}=1$, $\delta_{D}=\delta_{\Sigma}=0.5$.}
	\label{errtest1}
\end{figure}

Errors indicators $\mathcal{E}_{L^2}$, $\mathcal{E}_{H^1}$ are chosen for the 3D solution and $\widehat{\mathcal{E}}_{L^2}$ and $\widehat{\mathcal{E}}_{H^1}$ for the 1D problem, defined as:
\begin{align*}
&\mathcal{E}_{L^2}=\cfrac{||u_{ex}-U||_{L^2(\Omega)}}{||u_{ex}||_{L^2(\Omega)}} ,\qquad \mathcal{E}_{H^1}=\cfrac{||u_{ex}-U||_{H^1(\Omega)}}{||u_{ex}||_{H^1(\Omega)}},\\
&\widehat{\mathcal{E}}_{L^2}=\cfrac{||\hat{u}_{ex}-\hat{U}||_{L^2(\Lambda)}}{||\hat{u}_{ex}||_{L^2(\Lambda)}}, \qquad \widehat{\mathcal{E}}_{H^1}=\cfrac{||\hat{u}_{ex}-\hat{U}||_{H^1(\Lambda)}}{||\hat{u}_{ex}||_{H^1(\Lambda)}}.
\end{align*}
Figure~\ref{errtest1} displays the convergence trends for the above quantities against mesh refinement. Four meshes are considered, obtained by choosing $h=0.208,0.131,$ $0.083,0.052$, which correspond to $N=257,1026,4155,16545$ DOFs and $\hat{N}=15,29,57,88$ DOFs, respectively.  The parameters $\hat{\delta}_u=1$ and $\delta_D=\delta_\Sigma=0.5$ are fixed for all cases, hence a 3D mesh refinement induces a refinement of all the 1D meshes. Optimal convergence rates are obtained for the considered indicators in relation to the regularity of the exact solution, as reported in the picture. 

\begin{figure}
	\centering
	\includegraphics[width=0.45\linewidth]{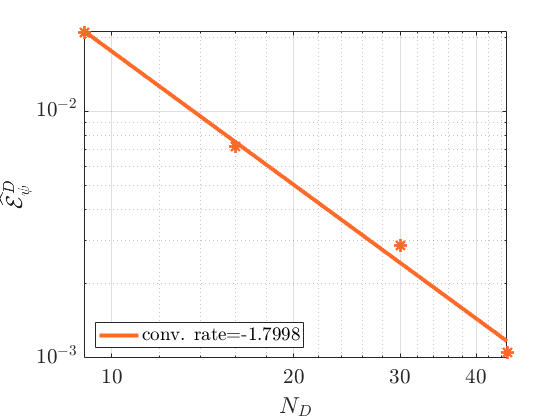}%
	\includegraphics[width=0.45\linewidth]{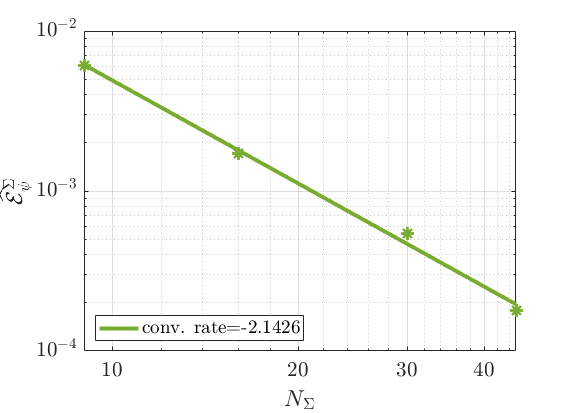}
	\caption{{TP1: trend of the error indicators for $\Psi_D$ and $\Psi_\Sigma$ under mesh refinement. On the left $\widehat{\mathcal{E}}_{\psi}^D$, on the right $\widehat{\mathcal{E}}_{\psi}^\Sigma$. Other parameters: $\hat{\delta}_{u}=1$, $\delta_{D}=\delta_{\Sigma}=0.5$.}}
	\label{errPhiPsi}
\end{figure}

Two additional error indicators are instead considered for the interface variables $\Psi_D$ and $\Psi_\Sigma$ as:
\begin{equation*}
\widehat{\mathcal{E}}_{\psi}^D=\cfrac{||\check{u}_{ex}-\Psi_D||_{L^2(\Lambda)}}{||\check{u}_{ex}||_{L^2(\Lambda)}}., \qquad \widehat{\mathcal{E}}_{\psi}^\Sigma=\cfrac{||\hat{u}_{ex}-\Psi_\Sigma||_{L^2(\Lambda)}}{||\hat{u}_{ex}||_{L^2(\Lambda)}}.
\end{equation*}
where $\extg\check{u}_{ex}=\trg u_{ex}=u_{ex}(\hat{R},z)$. The values of these errors indicators on the same meshes considered before, are reported in Figure~\ref{errPhiPsi}, on the left for $\widehat{\mathcal{E}}_{\psi}^D$ and, on the right, for $\widehat{\mathcal{E}}_{\psi}^\Sigma$.
Again optimal convergence curves are obtained, if it is considered that $\Psi_D$ is the approximation of a quantity in $H^{\frac12}(\Lambda)$ and $\Psi_\Sigma$ the approximation of a function in $H^1(\Lambda)$.

\begin{figure}
	\centering
	\includegraphics[width=0.45\linewidth]{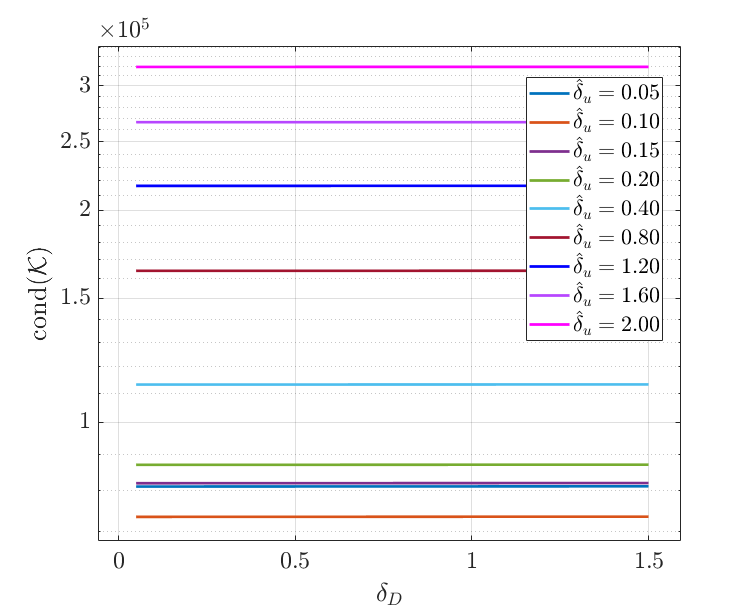}%
	\includegraphics[width=0.45\linewidth]{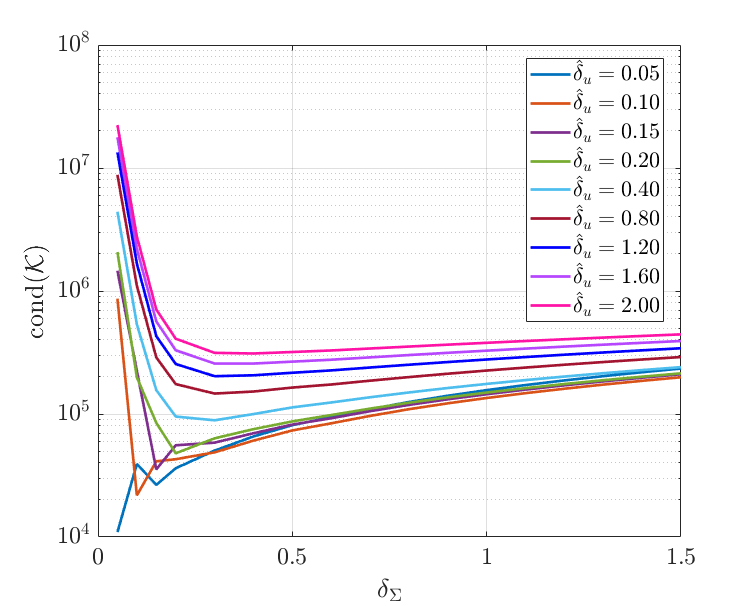}
	\caption{TP1: trend of the conditioning of the KKT system under the variation of the 1D mesh parameters. On the left variable $\delta_D$ and different values of $\hat{\delta} u$ while $\delta_{\Sigma}=0.5$. On the right variable  $\delta_{\Sigma}$ and $\delta_D=0.5$. In both cases $h=0.083$}
	\label{condkkt}
\end{figure}

\begin{figure}
	\centering
	\includegraphics[width=0.45\linewidth]{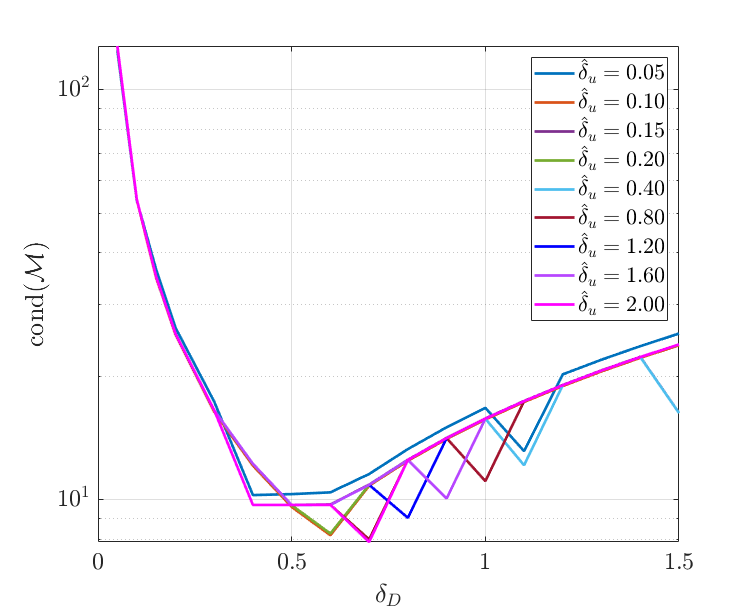}%
	\includegraphics[width=0.45\linewidth]{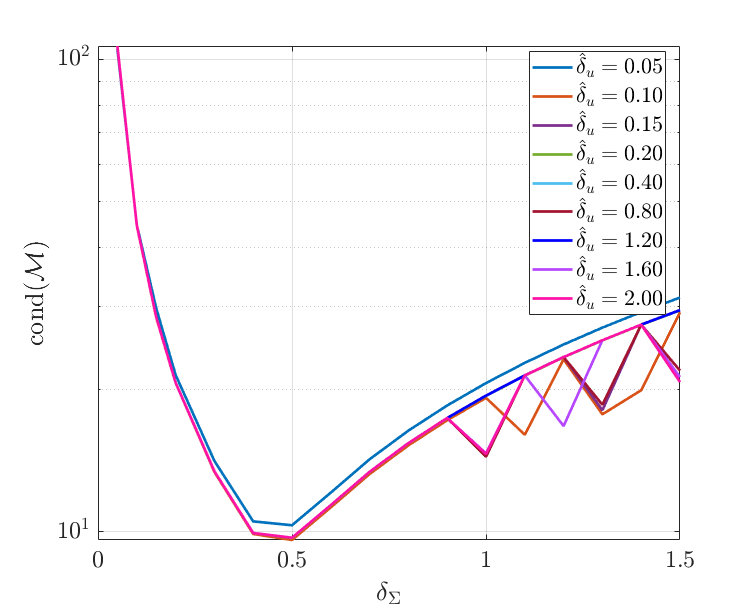}
	\caption{TP1: trend of the conditioning of system \eqref{systCg} under the variation of the 1D mesh parameters. On the left variable $\delta_D$ and different values of $\hat{\delta} u$ while $\delta_{\Sigma}=0.5$. On the right variable  $\delta_{\Sigma}$ and $\delta_D=0.5$. In both cases $h=0.083$}
	\label{condm}
\end{figure}

Figure~\ref{condkkt} shows the trend of the condition number of matrix $\mathcal{K}$, defined in \eqref{KKT}, under the variation of the 1D-mesh parameters. On the left the conditioning is plotted under the variation of $\delta_D$ and for different values of $\hat{\delta}_u$, while a constant $\delta_{\Sigma}=0.5$ is used; on the right $\delta_{\Sigma}$ varies instead, while $\delta_D=0.5$ is fixed, again for $\hat{\delta}_u$ ranging between 0.05 and 2. In both cases we can observe how, in general, at the increase $\hat{\delta}_u$, slightly higher conditioning values are registered,  with some exceptions for very small values of the parameter. In all cases, however, the impact of this parameter is quite marginal on the conditioning of the system. Looking at the left plot we can see that the value of $\delta_D$ has no impact on the conditioning.
Looking instead at Figure~\ref{condkkt} on the right, it can be noticed that $\delta_\Sigma$ has a larger impact on the conditioning, but only if very small values are used, and, at the same time, a value $\hat{\delta}_u>\delta_\Sigma$ is chosen; in these cases an increase of up to two orders of magnitude in the conditioning is observed. However, for $\delta_\Sigma>0.2$ the effect of $\delta_\Sigma$ on the conditioning becomes almost irrelevant, independently from the choice of the other parameters. The behaviour here observed is quite different from the one observed in Ref.~\cite{BGS3D1D2022}, where a larger effect of the mesh parameters on the conditioning was instead observed. Further, system \eqref{systCg} is known to be even better conditioned than the corresponding system \eqref{KKT}, \cite{Pestana2016}. For this simple example it is possible to explicitly compute matrix $\mathcal{M}$. Its conditioning is plotted in Figure~\ref{condm}, against variations of the 1D mesh parameters. Trends similar to the ones of Figures~\ref{condkkt} are observed and the behaviour is almost the same for all the values of $\hat{\delta}_u$ considered, but the conditioning of $\mathcal{M}$ is, in general, between 4 and 5 orders of magnitude smaller than the one of $\mathcal{K}$. This is expected to have a positive impact on the number of iterations of Algorithm~\ref{grad_congAlgo}, with few iterations required to reach the prescribed tolerance even without the use of a preconditioner. The analysis of the performances of the proposed iterative solver is deferred to the last example here proposed, in which a more complex setting is considered.


\subsection{Test Problem 2 (TP2)}
\begin{figure}
	\centering
	\includegraphics[width=0.55\linewidth]{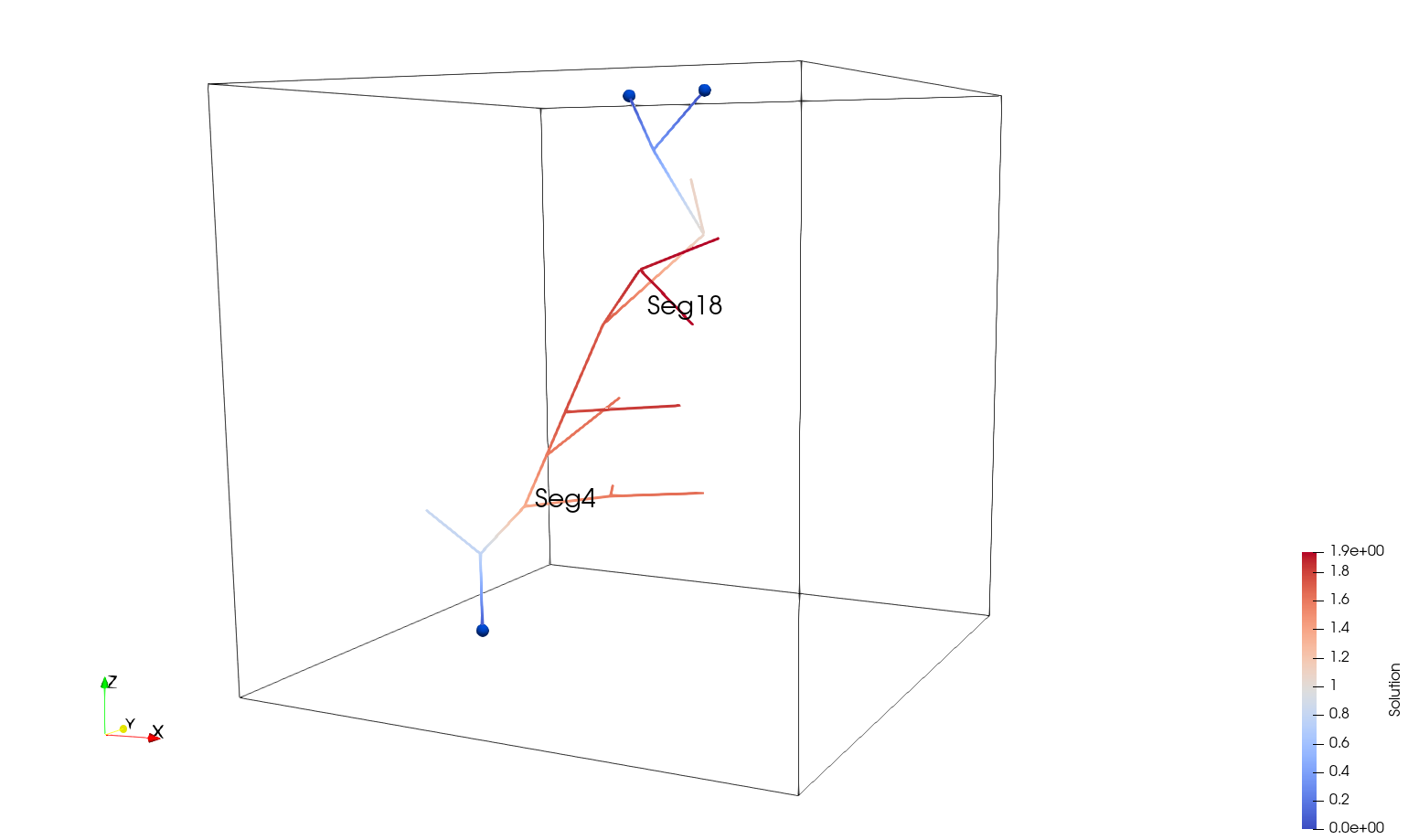}
	\caption{TP2: Solution obtained on the centerlines of the inclusions for $h=0.083$, $\hat{\delta}_u=1$ and $\delta_{D}=\delta_{\Sigma}=0.5$. The point marked in blue identify homogeneous Dirichlet boundary conditions.}
	\label{config19seg}
\end{figure}

For this numerical example we consider a set of $19$ inclusions of radius $\check{R}=10^{-2}$, whose centerlines intersect in $9$ points. The resulting network is embedded in the same cube of edge $l=2$ considered in Test Problem 1. We impose homogeneous Dirichlet boundary conditions on all the faces of the cube and at the dead ends of the network which intersect the cube at its top and bottom faces, as shown in Figure~\ref{config19seg}. Homogeneous Neumann boundary conditions are imposed at  segment endpoints lying inside the cube. For what concerns problem coefficients we consider ${K}=1$, $f=0$ and ${\tilde{K}}_i=100$, $\overline{\overline{g}}_i=100$, $\beta_i=5e-2$, $\forall i=1,...,19$.
\begin{figure}
	\centering
	\includegraphics[width=0.5\linewidth]{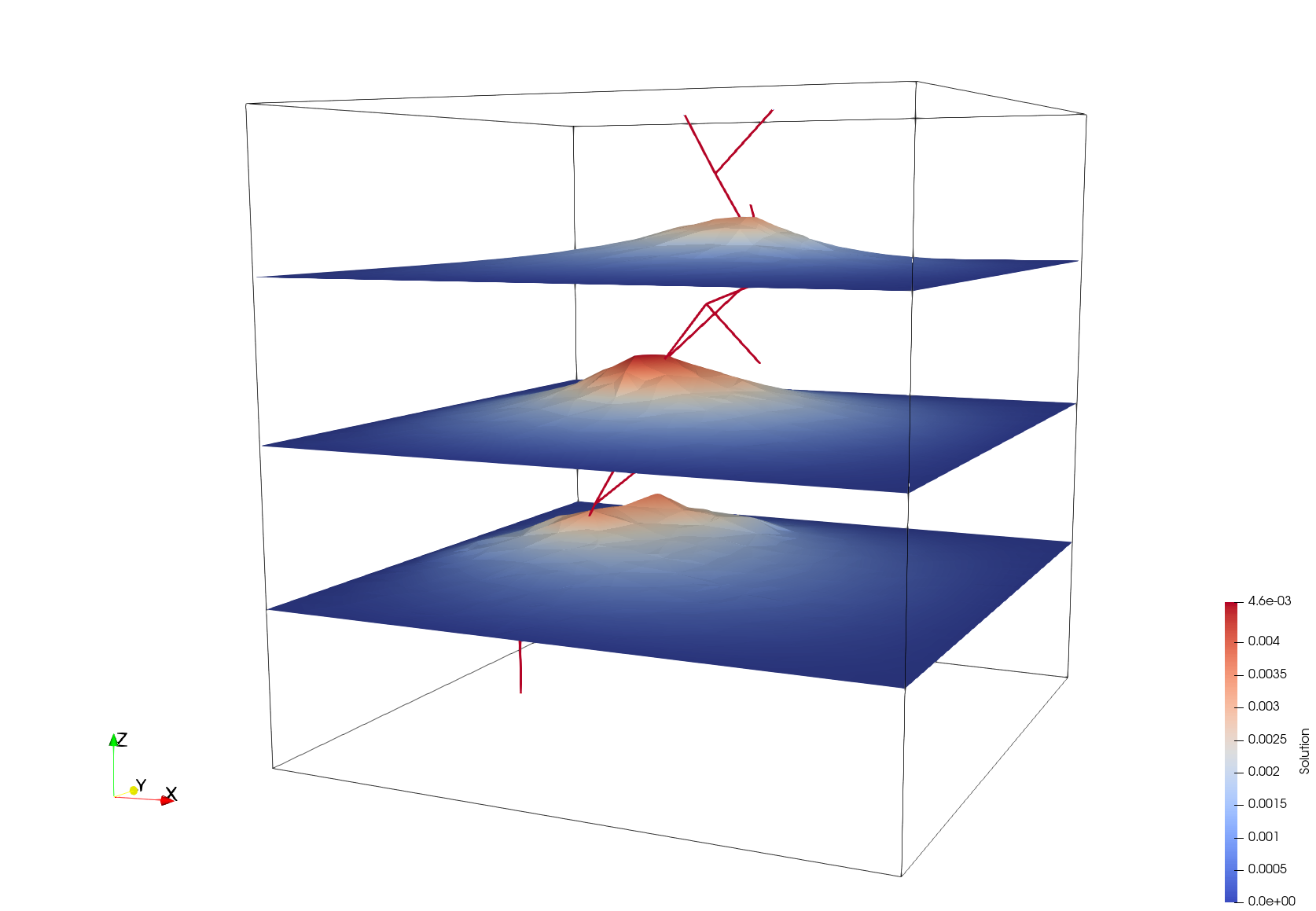}%
	\includegraphics[width=0.5\linewidth]{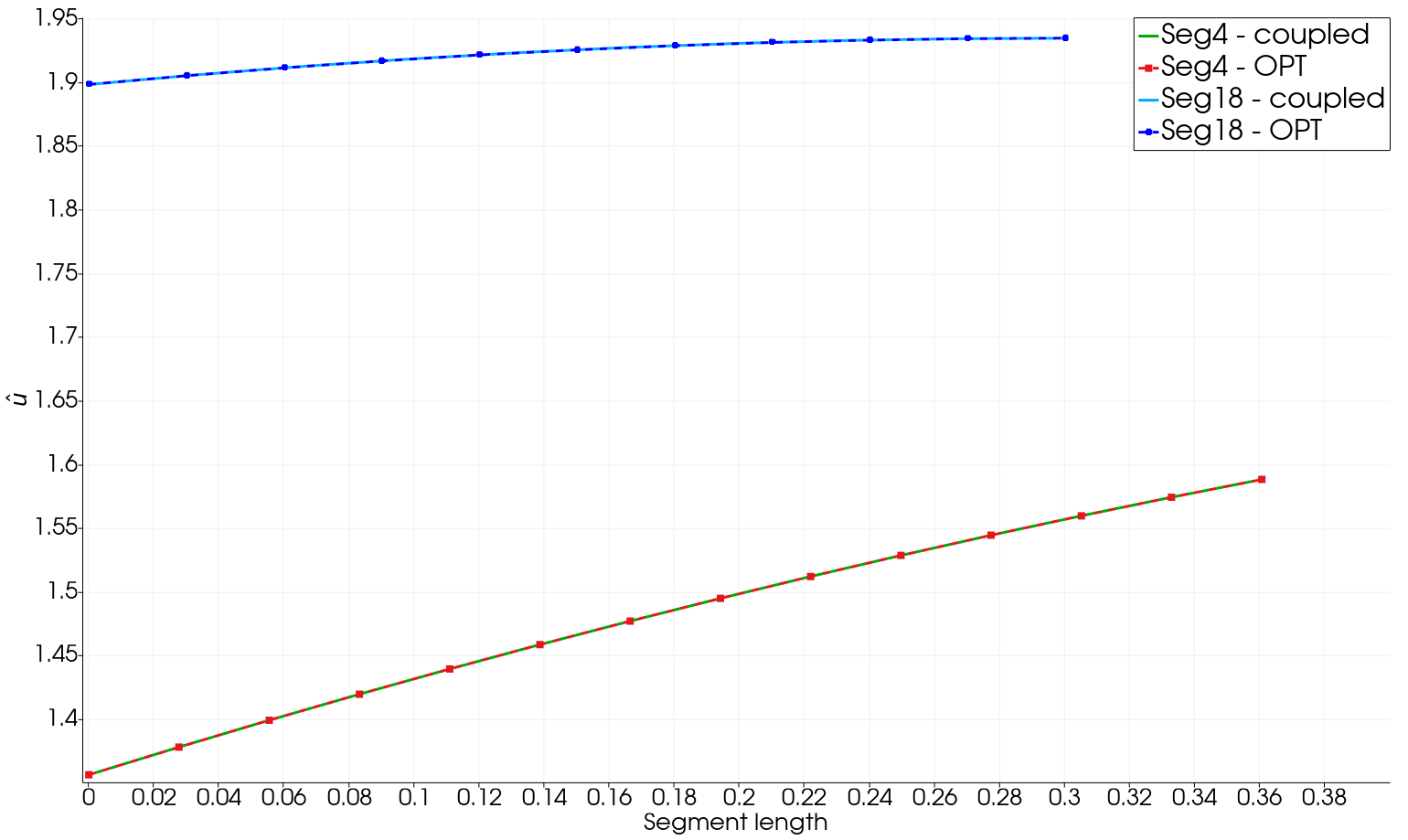}
	\caption{TP2: On the left: solutions with OPT and coupled methods obtained inside the cube on three different planes parallel to the $x-y$ plane and located at $z=-0.5$, $z=0$ and $z=0.5$. Solution amplified by a factor $50$ with respect to domain size; on the right: comparison of the solution on two selected segments with the OPT and the coupled methods.}
	\label{solslice}
\end{figure}

The problem is solved both with the method proposed in this article and with a different approach in which no auxiliary variables are introduced. The following 3D-1D coupled problem is derived from \eqref{coup_eq1}-\eqref{coup_eq2}: \textit{find} $(u, \hat{u}) \in V_D \times \hat{V}$ \textit{such that}:
\begin{align*}
&({{K}}\nabla u, \nabla v)_{L^2(D)}+\left( |\Gamma|\beta\check{u},\check{v}\right) _{L^2(\Lambda)}-\left( |\Gamma|\beta\hat{u},\check{v} \right) _{L^2(\Lambda)}=(f,v)_{L^2(D)}  \\[-0.5em]  &\hspace{8.5cm}\forall v \in {V_D}, \check{v} \in \hat{V}: \gamma_{_\Gamma} v=\mathcal{E}_{_\Gamma}\check{v} \nonumber \\
&\Big( {\tilde{K}}|\Sigma|\frac{d\hat{u}}{ds},\frac{d\hat{v}}{ds}\Big)_{{L^2(\Lambda)}}+\left( |\Gamma|\beta\hat{u},\hat{v} \right) _{L^2(\Lambda)}-\left(  |\Gamma|\beta\check{u},\hat{v}\right) _{L^2(\Lambda)}=(|\Sigma|\overline{\overline{g}},\hat{v})_{{L^2(\Lambda)}} \qquad \forall \hat{v} \in \hat{V},
\end{align*}
which, after discretization yields the following global system
\begin{displaymath}
\begin{bmatrix}
\bm{A} &-\bm{B}\\
-\bm{B}^T & \bm{\hat{A}}
\end{bmatrix}\begin{bmatrix}
U \\ \hat{U}
\end{bmatrix}=\begin{bmatrix}
f \\ g
\end{bmatrix}
\end{displaymath}
where the nomenclature is the same of Section~\ref{Discrete} and the new matrix $\bm{B}$ is defined as
\begin{displaymath}
\bm{B}=\left[ \bm{B_1}, \bm{B_2},...,\bm{B_{\I}}\right]  \in \mathbb{R}^{N\times \hat{N}}
\end{displaymath}
with 
\begin{displaymath}
\bm{B_i} \in \mathbb{R}^{N\times \hat{N}_i} \text{ s.t. } (B_i)_{kl}=\int_{\Lambda_i}{{\varphi_{k}}_{|_{\Lambda_i}}\hat\varphi_{i,l}}~ds
.\end{displaymath}
We will refer to this method as ``coupled'' and we will use it as a comparison term for our approach, which, instead, will be labelled as OPT (optimization based). The results shown in the following are obtained by considering a 3D mesh, non conforming to the inclusions, with $h=0.083$ and $N=3320$ and a 1D mesh with $\hat{\delta}_{u}=1$, corresponding to $\hat{N}=234$ DOFs. Parameters  $\delta_{D}=\delta_{\Sigma}=0.5$ are used in the OPT approach.
Figure~\ref{config19seg} reports the solutions obtained with both approaches on the  network of segments, while Figure \ref{solslice} on the right proposes a comparison of the solutions on two selected segments (marked in Figure~\ref{config19seg}), showing an almost perfect agreement. Figure~\ref{solslice} on the left, instead, reports the solution $U$ obtained inside the cube on three different planes, located at $z=-0.5$, $z=0$ and $z=0.5$ and all parallel to the $x-y$ plane. Even in this Figure the solutions obtained with both the approaches are shown and appear to be almost perfectly overlapped.

\subsection{Conjugate gradient test (CGtest)}
\begin{figure}
	\centering
	\includegraphics[width=0.8\linewidth]{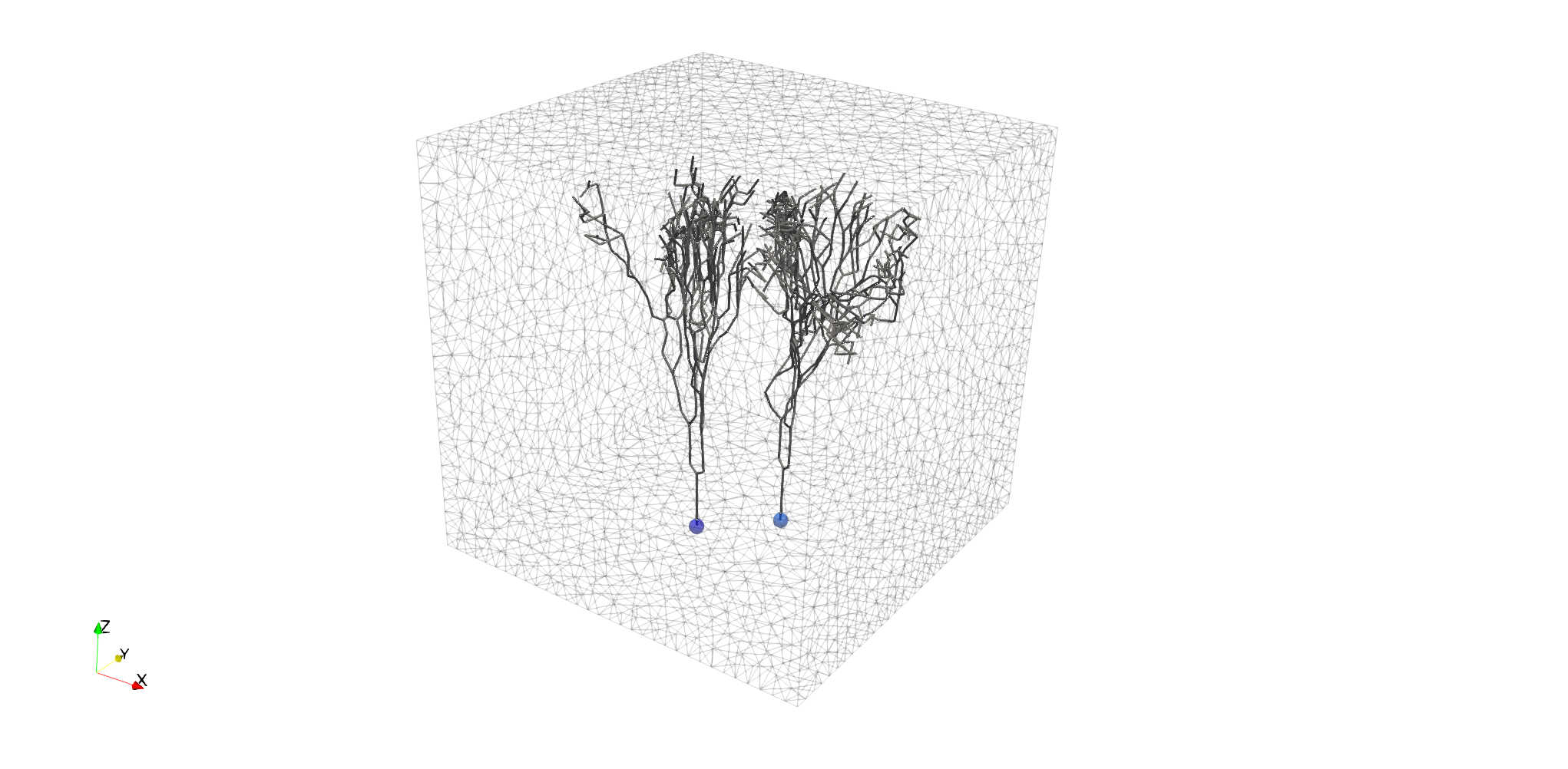}
	\caption{CGtest: Representation of the problem geometry. The blue spheres highlight the Dirichlet boundary conditions at the inlets of the two networks.}
	\label{MIgeometry}
\end{figure}

\begin{figure}
	\centering
	\includegraphics[width=0.6\linewidth]{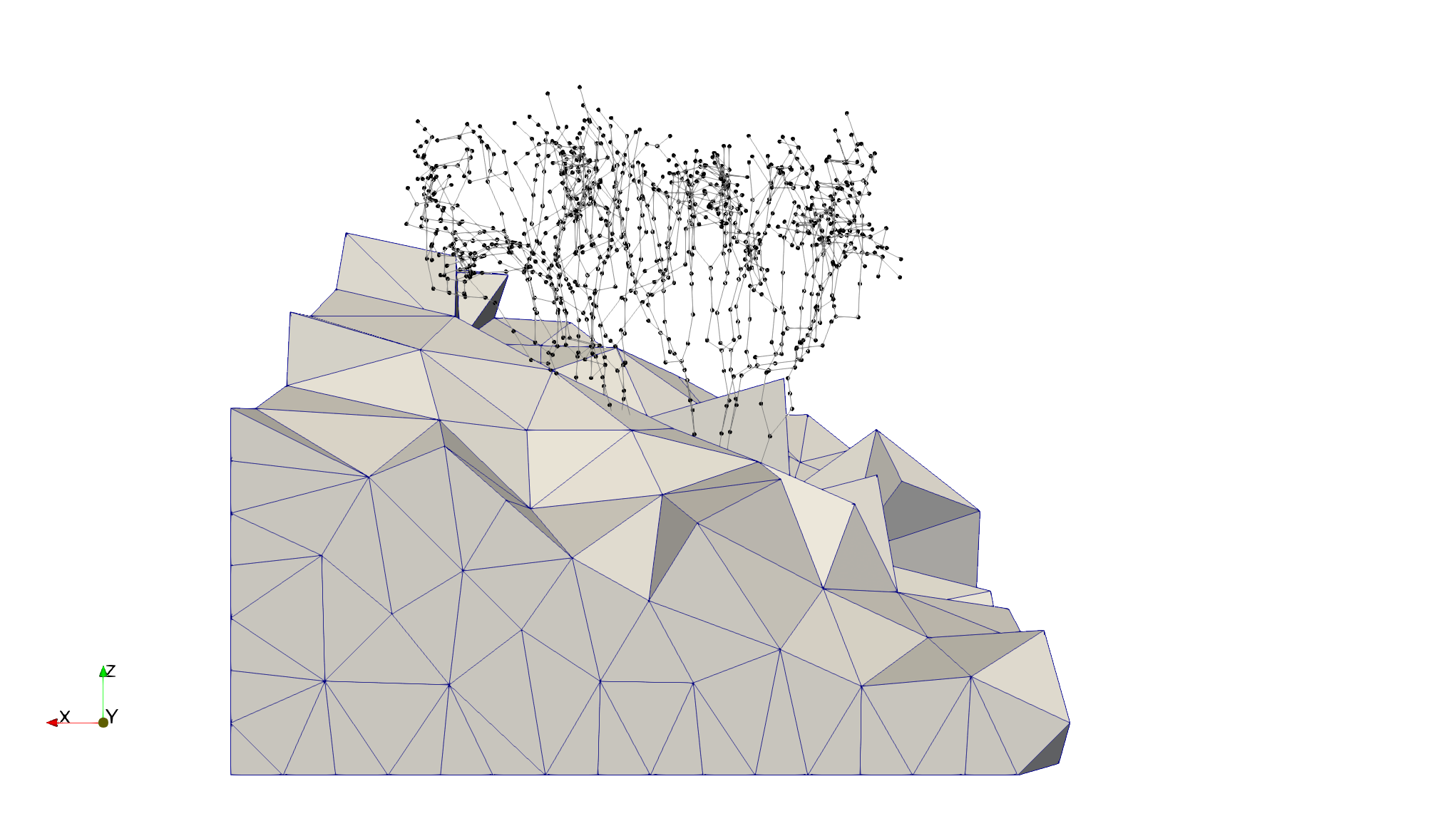}\hspace{-3.1cm}%
	\includegraphics[width=0.6\linewidth]{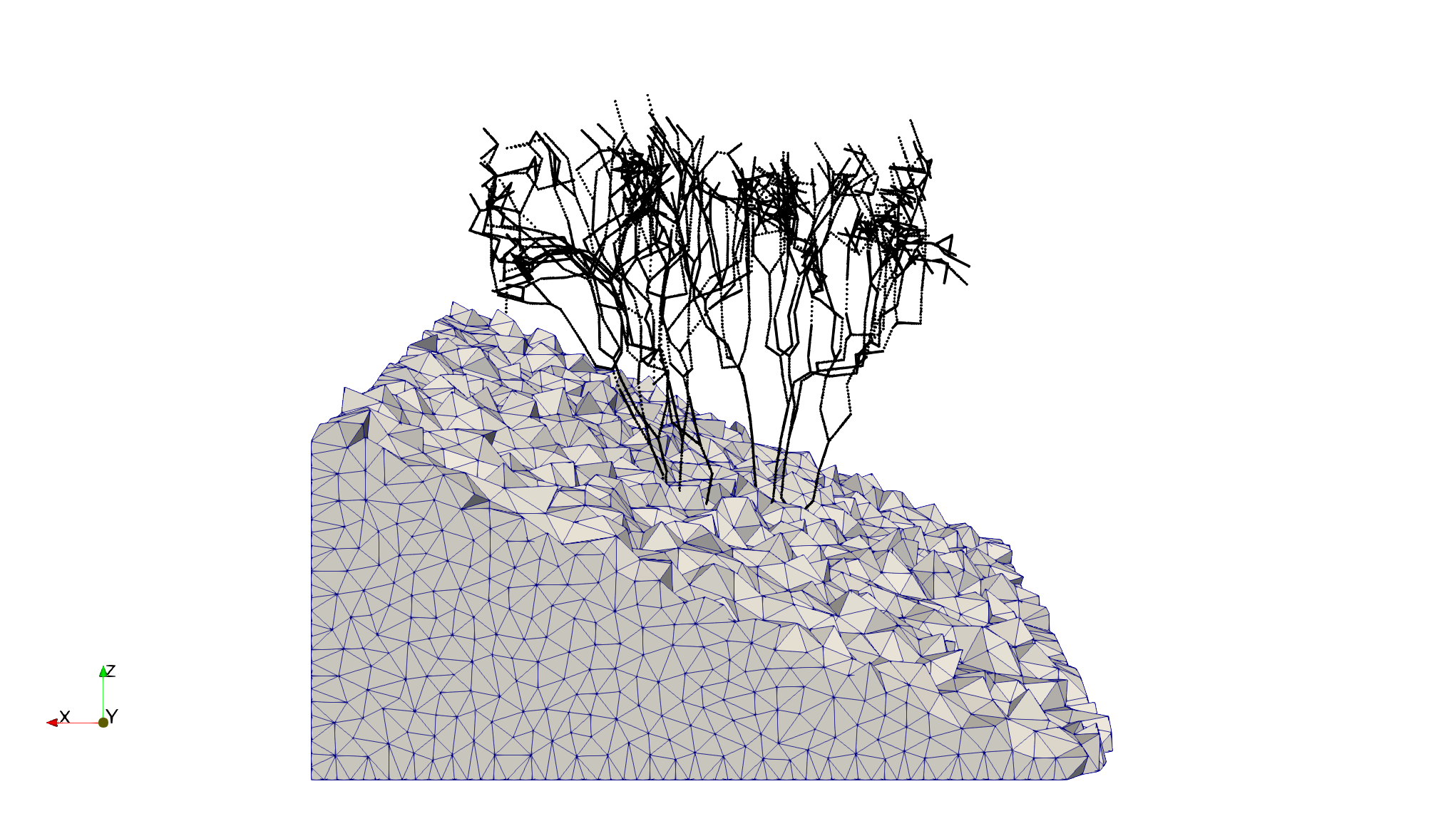}
	\caption{CGtest: comparison between meshes. On the left, $h=0.208$, $\delta_{D}=\delta_{\Sigma}=0.5$; on the right, $h=0.052$, $\delta_{D}=\delta_{\Sigma}=0.5$. }
	\label{mesh}
\end{figure}

\begin{figure}
	\centering
	\includegraphics[width=0.5\linewidth]{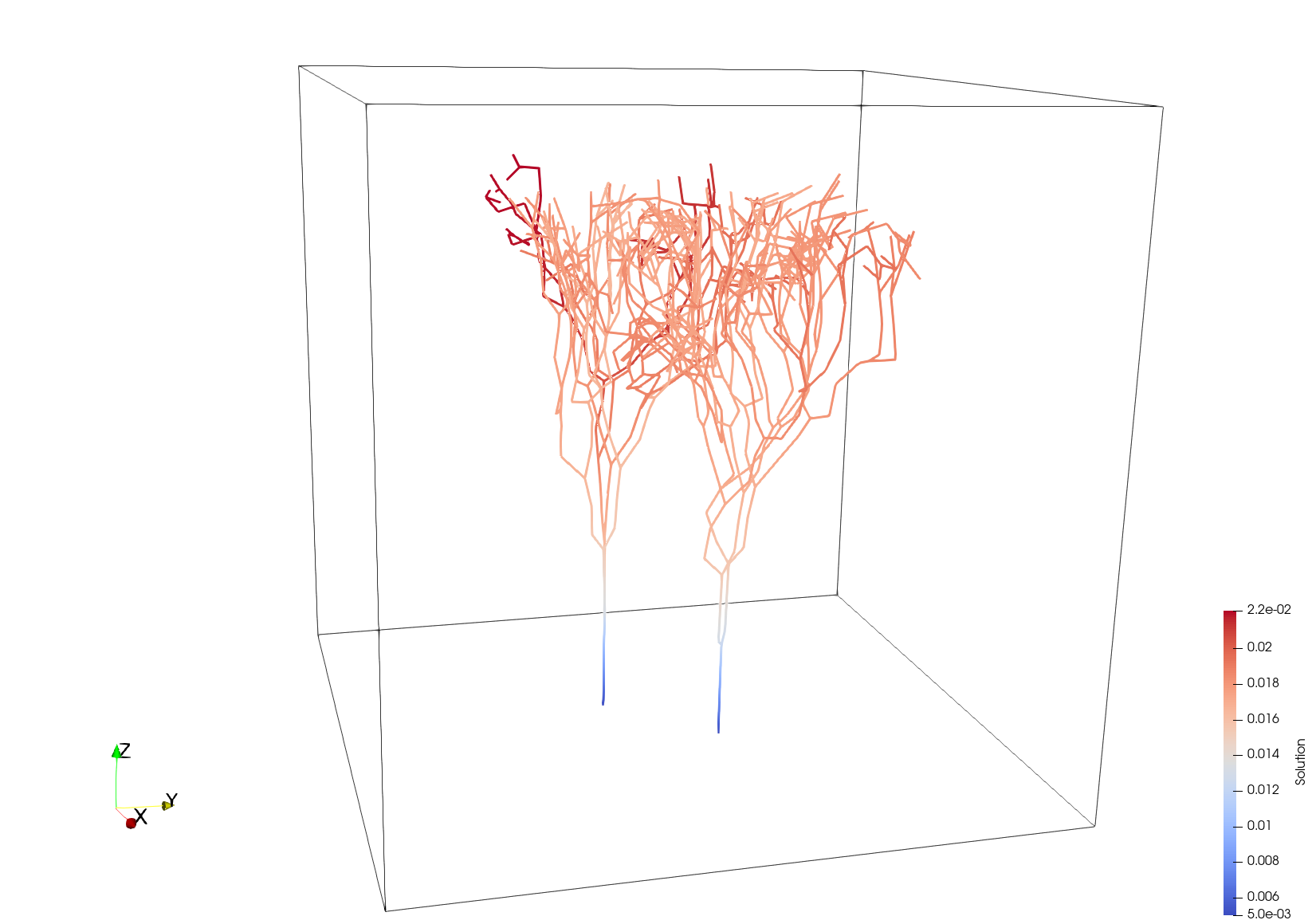}%
	\includegraphics[width=0.5\linewidth]{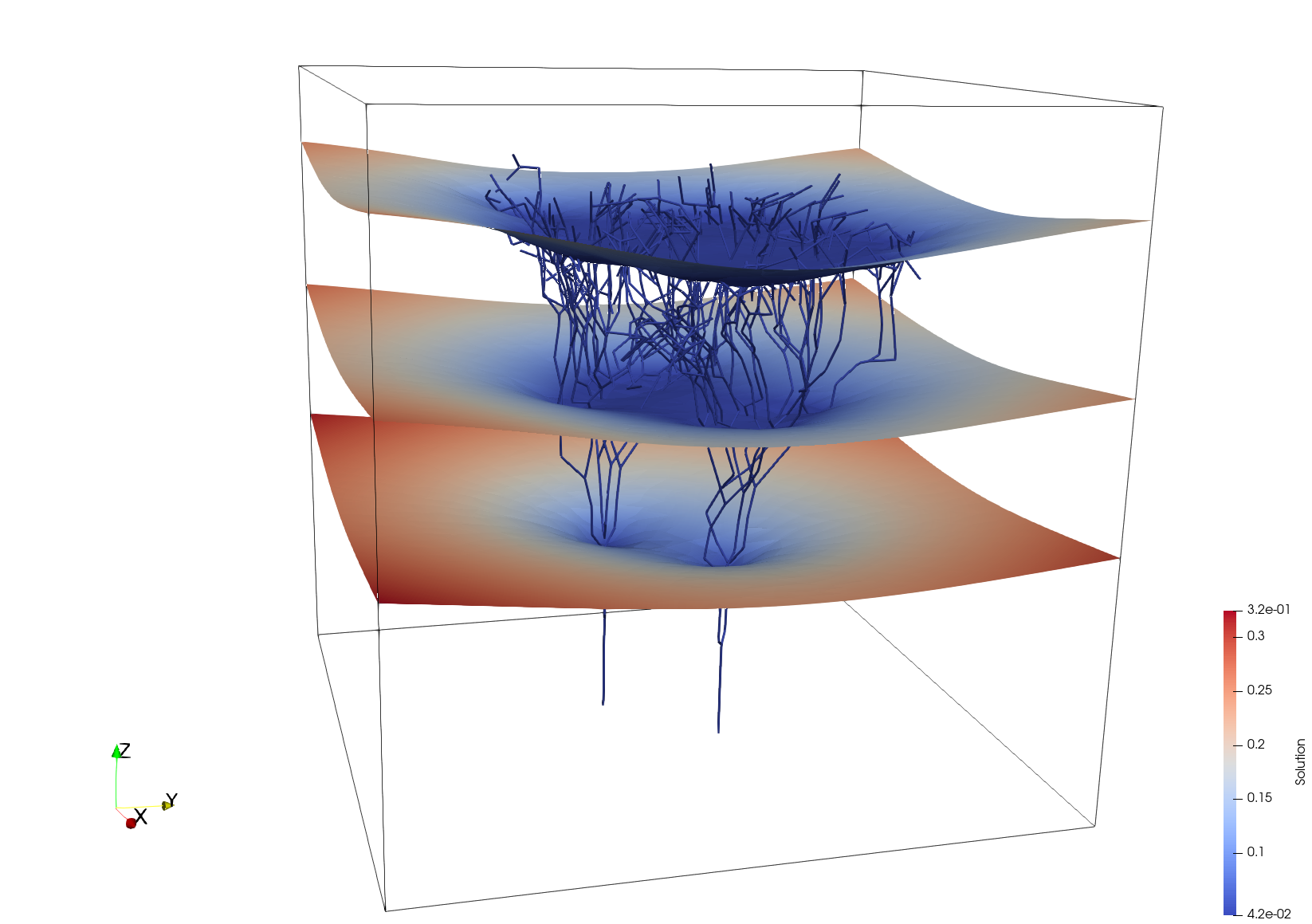}
	\caption{CGtest: on the left, solution obtained inside the inclusions; on the right, solution obtained in the cube on three different planes perpendicular to the $z$-axis, namely $z=-0.5$, $z=0$ and $z=0.5$. Parameters: $h=0.052$, $\delta_{D}=\delta_{\Sigma}=2$, $\hat{\delta}_u=1$.}
	\label{sol}
\end{figure}

We now consider a more complex numerical example, characterized by the presence of multiple intersecting inclusions. The setting of this example might be considered as realistic of a living tissue with a network of vessels. The purpose of the present example is to test the performances of the proposed resolution strategy and preconditioner in a realistic setting. In particular we consider the domain of Figure~\ref{MIgeometry}, where $873$ segments organized into two connected clusters are immersed in a cubic domain $\Omega$ of edge $l=2$, as considered in the previous examples. On the faces of the cube we consider Neumann boundary conditions, namely $K\nabla u\cdot \bm{n}=2\cdot 10^{-5}, $ with $\bm{n}$ denoting, in this case, the outward pointing unit normal vector to $\partial \Omega$. At the inlets of the two networks, i.e. the segment endpoints lying on the face $z=-1$, we impose Dirichlet boundary conditions equal to $5\cdot 10^{-3}$, and at all other segment endpoints we consider homogeneous Neumann conditions. Problem data are as follows: 
\begin{equation*}
K=2\cdot 10^{-4}, \qquad f(x,y,z)=0, \qquad\tilde{K}=3\cdot 10^{1},\qquad \overline{\overline{g}}=0\qquad \beta=1\cdot 10^{-2}.
\end{equation*}
The cluster of segments of this example is composed by many segments of small length, compared to domain size. Thus, a set of simulations is performed by varying both parameter $h$ of the 3D mesh and parameters $\delta_{D}$ and $\delta_{\Sigma}$ of the 1D meshes, in order to control the meshsize independently, being, instead, $\hat{\delta}_u=1$ fixed. The values of the parameters used in the simulations are reported in the first three columns of Table~\ref{itcg}, whereas the fourth column reports the number of DOFs of the variables $\Psi_D$ and $\Psi_\Sigma$, corresponding to the size of the system \eqref{systCg}. The coarsest and the finest mesh combination considered are shown in Figure~\ref{mesh}, whereas the solution on the finest mesh is reported in Figure~\ref{sol}, on the left for the 1D network and on the right for the 3D solution on three planes orthogonal to the $z$-axis. The remaining columns of Table~\ref{itcg} report the number of iterations required by Algorithm~\ref{grad_congAlgo} to solve system \eqref{systCg} up to a relative residual of $10^{-6}$ or $10^{-9}$, without using a preconditioner (columns $\bm{CG_{it}}^{(10^{-6})}$ and $\bm{CG_{it}}^{(10^{-9})}$) or using preconditioner \eqref{prec} (columns $\bm{PCG_{it}}^{(10^{-6})}$ and $\bm{PCG_{it}}^{(10^{-9})}$).
\begin{table}\label{itcg}
	\caption{CGtest: mesh parameters, DOFs and corresponding CG iterations. In brackets value of the relative residual which defines the stopping criterion. For all the considered cases $\hat{\delta}_u=1$.}
	{
		\begin{tabular}{@{}cccccccc@{}} \hline
			$h$ & $\delta_{D}$&$\delta_{\Sigma}$&  $N_D+N_\Sigma$ &$\bm{CG_{it}}^{(10^{-6})}$  &$\bm{CG_{it}}^{(10^{-9})}$  & $\bm{PCG_{it}}^{(10^{-6})}$ & $\bm{PCG_{it}}^{(10^{-9})}$\\
			\hline
			$0.208$&$0.5$& $0.5$ & $3650$ &$39$ &$57$ & $33$ &$43$\\
			$0.131$&$1.0$& $1.0$ & $6344$ &$48$ &$67$
			&$35$ &$46$\\
			$0.083$&$1.5$& $1.5$ & $12428$ &$47$ &$68$ &$36$ &$48$\\
			$0.052$&$2.0$& $2.0$ & $21256$ &$44$ &$61$ &$37$ &$49$\\\hline
	\end{tabular}}
\end{table}

It can be seen that, in all cases, the number of iterations is small compared to the number of unknowns, and it only marginally grows as the stopping tolerance is reduced. The use of the preconditioner allows to further reduce the number of iteration, the obtained reduction ranging between $15\%$ and $30\%$. It is to be highlighted that the proposed preconditioner can be obtained and applied at a very low computational cost, as it only involves the resolution of local 1D problems and can be computed and used in parallel. The effectiveness of the proposed resolution approach reflects the good conditioning of the obtained system, as it was pointed out in Test Problem 1.


\section{Conclusions}\label{Concl}
A PDE constrained formulation for 3D-1D coupled problems with discontinuous solution at the interfaces has been derived and proposed. The approach is based on the introduction of unknown interface variables to decouple the 3D and 1D problems and on the minimization of a cost functional to enforce interface conditions. The problem is discretized resorting to standard finite elements on non conforming meshes independently set on each subdomain. Well posedness results for the discrete problem are obtained independently of the choice of the mesh parameters of the various domains. The proposed test on a problem with known analytical solution shows that optimal convergence trends of the error are obtained for both the 3D and 1D solution. Also the linear system corresponding to the application of the method appears to be well conditioned for a wide range of choices of the mesh parameters. The examples on more complex domains reveal the applicability of the method to realistic configurations and also the good performances of the proposed gradient-based solver.

\section*{Acknowledgements}
This work is supported by the MIUR project ``Dipartimenti di Eccellenza 2018-2022'' (CUP E11G18000350001), PRIN project ``Virtual Element Methods: Analysis and Applications'' (201744KLJL\_004) and by INdAM-GNCS. Computational resources are partially supported by SmartData@polito.

\printbibliography

\end{document}